\documentclass[10pt]{amsart}
\pdfoutput=1
\usepackage{amsmath}
\usepackage{amssymb}
\usepackage{amsthm}
\usepackage{bm}
\usepackage{accents}
\usepackage{mathtools}
\usepackage{tikz}
\usetikzlibrary{calc}
\usetikzlibrary{decorations.pathmorphing,shapes}
\usetikzlibrary{automata,positioning,quotes}
\usepackage{tikz-cd}
\usepackage{forest}
\usepackage{braket} 
\usepackage{listings}
\usepackage{mdframed}
\usepackage{verbatim}
\usepackage{physics2}
\usephysicsmodule{ab,ab.legacy,diagmat,xmat}
\usepackage{derivative}
\usepackage{fixdif}
\usepackage{stmaryrd}
\usepackage[mathscr]{eucal}
\usepackage{stackengine} 

\usepackage{microtype}

\usepackage{algorithmicx}
\usepackage{algpseudocode}
\usepackage{algorithm}

\usepackage[english]{babel} 
\usepackage[backend=biber,style=alphabetic,maxalphanames=4,maxnames=5,hyperref]{biblatex}
\usepackage[bookmarks, colorlinks, breaklinks]{hyperref} 
\hypersetup{linkcolor=blue,citecolor=magenta,filecolor=black,urlcolor=blue}
\usepackage{cleveref}
\usepackage{xurl}
\usepackage{graphicx}
\graphicspath{{./}}
\crefname{equation}{}{}

\usepackage{float}
\usepackage{booktabs}
\usepackage[shortlabels]{enumitem}
\setitemize{noitemsep}
\usepackage{csquotes}
\newlist{mydescription}{description}{1}
\setlist[mydescription]{style=nextline,
                        font=\bfseries,
                        labelindent=1cm, 
                        leftmargin =2cm,
                        rightmargin=1cm,
                        topsep     =1ex
                       }

\usepackage{lipsum}

\DeclarePairedDelimiter{\floor}{\lfloor}{\rfloor}
\DeclarePairedDelimiter{\ceil}{\lceil}{\rceil}

\newtheorem{thm}{Theorem}[section]
\newtheorem{cor}[thm]{Corollary}
\newtheorem{prop}[thm]{Proposition}
\newtheorem{lem}[thm]{Lemma}

\newtheorem{fthm}[thm]{``Theorem''}

\theoremstyle{definition}
\newtheorem{defn}[thm]{Definition}
\newtheorem{defns}[thm]{Definitions}

\newtheorem{exm}[thm]{Example}

\newtheorem{notn}[thm]{Notation}

\newtheorem{conv}[thm]{Convention}

\theoremstyle{remark}
\newtheorem{rmk}[thm]{Remark}
\newtheorem{rmks}[thm]{Remarks}
\newtheorem{warn}[thm]{Warning}

\theoremstyle{plain}
\newtheorem*{thm*}{Theorem}
\newtheorem*{prop*}{Proposition}
\newtheorem*{lem*}{Lemma}
\newtheorem*{cor*}{Corollary}
\newtheorem*{conj*}{Conjecture}

\theoremstyle{definition}
\newtheorem*{defn*}{Definition}
\newtheorem*{exer*}{Exercise}
\newtheorem*{defns*}{Definitions}
\newtheorem*{con*}{Construction}
\newtheorem*{exm*}{Example}
\newtheorem*{exms*}{Examples}
\newtheorem*{notn*}{Notation}
\newtheorem*{notns*}{Notations}
\newtheorem*{addm*}{Addendum}

\theoremstyle{remark}
\newtheorem*{rmk*}{Remark}

\newcommand{\A}{\mathbb{A}}
\newcommand{\G}{\mathbb{G}}

\newcommand{\R}{\mathbb{R}}
\newcommand{\C}{\mathbb{C}}
\newcommand{\Z}{\mathbb{Z}}
\newcommand{\Q}{\mathbb{Q}}
\newcommand{\F}{\mathbb{F}}
\newcommand{\W}{\mathcal{W}}
\newcommand{\bd}{\mathbf{d}}
\newcommand{\ba}{\mathbf{a}}

\renewcommand{\P}{\mathbb{P}}
\newcommand{\E}{\mathbb{E}}

\newcommand{\ep}{\varepsilon}

\newcommand{\mc}[1]{\mathcal{#1}}
\newcommand{\msc}[1]{\mathscr{#1}}
\newcommand{\T}{\mathbf{T}}
\newcommand{\mf}[1]{\mathfrak{#1}}
\newcommand{\mbf}[1]{\mathbf{#1}}
\newcommand{\mr}[1]{\mathrm{#1}}
\newcommand{\on}[1]{\operatorname{#1}}

\newcommand{\ol}[1]{\overline{#1}}

\newcommand{\wt}[1]{\widetilde{#1}}

\newcommand{\1}{\mathbf{1}}

\newcommand{\vir}{\mr{vir}}

\newcommand{\pt}{\mr{pt}}
\newcommand{\loc}{\mr{loc}}
\newcommand{\tw}{\mr{tw}}
\renewcommand{\top}{\mr{top}}
\newcommand{\ps}[1]{\llbracket #1 \rrbracket}

\DeclareMathOperator{\End}{End}

\DeclareMathOperator{\Aut}{Aut}

\DeclareMathOperator{\ev}{ev}
\DeclareMathOperator{\PF}{\msc{PF}}
\DeclareMathOperator{\st}{st}
\DeclareMathOperator{\Res}{Res}
\DeclareMathOperator{\cl}{cl}
\DeclareMathOperator{\Cont}{Cont}

\numberwithin{equation}{section}
\title[Polynomiality of higher genus GW]{Higher genus Gromov-Witten theory of one-parameter Calabi-Yau threefolds I: Polynomiality}
\author{Patrick Lei}
\date{\today}

\begin{document}
    
\begin{abstract}
    We prove the finite generation conjecture of~\cite{yy04} for the Gromov-Witten potentials of the Calabi-Yau hypersurfaces $Z_6 \subset \P(1,1,1,1,2)$, $Z_8 \subset \P(1,1,1,1,4)$, and $Z_{10} \subset \P(1,1,1,2,5)$ using the theory of MSP fields. In addition, a formula is given for the genus one Gromov-Witten potentials of these targets.
\end{abstract}

\maketitle

\tableofcontents

\section{Introduction}
\label{sec:introduction}

\subsection{Historical overview}%
\label{sub:Historical overview}

Mathematical interest in Gromov-Witten theory was sparked by the predictions of~\cite{cdgp}, which predicted the number of genus $0$ curves of all degrees (precisely, the genus $0$ Gromov-Witten potential $F_0$) in the quintic threefold using mirror symmetry, which is a duality between two superconformal field theories known as the $A$-model and the $B$-model.

Later, Gromov-Witten theory, which computes $A$-model invariants, was developed in symplectic topology by Ruan-Tian~\cite{rtqcoh} and in algebraic geometry by Behrend-Fantechi and Li-Tian~\cite{intrinsicnc,ltgwfoundation} and the mirror formula for $F_0$ of semi-positive complete intersections in projective toric varieties was proved by Givental and Lian-Liu-Yau~\cite{eqgwinv,giventalmirror,lly,lly2}. A key insight was to relate the Gromov-Witten theory of the complete intersection $Z \subseteq X_{\Sigma}$ defined by the vanishing of a section of a vector bundle $E$ with the so-called twisted GW theory of the pair $(X_{\Sigma}, E)$, which is known as the Quantum Lefschetz principle.

In higher genus, Bershadsky-Cecotti-Ooguri-Vafa~\cite{bcov} studied the $B$-model and deduced a \textit{Feynman rule}, which computes $F_g$ from the lower genus $F_{h<g}$ and a finite ambiguity. Later, Yamaguchi-Yau~\cite{yy04} studied the consequences of the BCOV conjecture in the $A$-model and deduced that
\begin{enumerate}
    \item All of the Gromov-Witten potentials $F_g$ for $Z_5 \subset \P^4$, $Z_6 \subset \P(1,1,1,1,2)$, $Z_8 \subset \P(1,1,1,1,4)$, and $Z_{10} \subset \P(1,1,1,2,5)$ lie in a ring $\msc{R}$ of polynomials in five power series in the K\"ahler parameter;
    \item The $F_g$ satisfy a set of partial differential equations called the Yamaguchi-Yau equations.\footnote{Often, mathematicians call these equations the Holomorphic Anomaly Equations. This is technically inaccurate because in physics, the $B$-model partition function takes the form $F_g^B(t, \ol{t})$ and has an antiholomorphic part. The physical Holomorphic Anomaly Equations govern the antiholomorphic part of $F_g^B(t, \ol{t})$.}
\end{enumerate}
Huang-Klemm-Quackenbush~\cite{hkq} used these results and new techniques to compute the Gromov-Witten invariants of the quintic for $g \leq 51$. Unfortunately, the mathematical formulation of BCOV theory is very difficult and has only been completed for the elliptic curve~\cite{mathbcov1,mathbcov2}.

Unfortunately, the Quantum Lefschetz principle fails in higher genus $g > 0$. Based on the idea of Li-Vakil-Zinger~\cite{desingularizationgenusone,g1cilizinger} to desingularize the ``main component'' of the moduli space of stable maps to $\P^4$, Zinger~\cite{reducedgenus1} proved the BCOV conjecture for the quintic in genus $1$. This result was later reproved by Kim-Lho~\cite{kimlho} using quasimap wall-crossing~\cite{qmapwcms}, which is another way to restore the Quantum Lefschetz principle in genus $1$.

A better way to restore the Quantum Lefschetz principle is to introduce the so-called $p$-fields, following Chang-Li~\cite{pfieldschangli}. They consider stable maps to $\P^4$, but with an additional section 
\[p \in H^0(\msc{C}, \msc{O}(-5) \otimes \omega_{\msc{C}}^{\log}). \]
The upshot is that the virtual cycle of Chang-Li's moduli space coincides with the virtual cycle of the moduli space of stable maps to the quintic up to an explicit sign. Using a compactification of this moduli space, Guo-Janda-Ruan~\cite{genus2logglsm, bcovlogglsm} proved first the BCOV conjecture for the quintic in genus $2$ and then in all genus, pending the proof of a virtual localization formula in log geometry.

The $p$-field construction realizes the Gromov-Witten theory of the quintic (and other toric complete intersections) as a mathematical version of a gauged linear sigma model (GLSM). GLSMs were introduced by Witten~\cite{mirrorandglsm} in physics and by Fan-Jarvis-Ruan~\cite{glsm} in mathematics. As a theory based on geometric invariant theory, a GLSM may have several phases corresponding to different GIT stability chambers. Physically, these correspond to different points on the so-called K\"ahler moduli space, which controls the $A$-model. The two phases of the GLSM for the quintic threefold are known as the \textit{geometric phase} and the \textit{Landau-Ginzburg phase}, which is a theory of the singularity $( [\C^5/\mu_5], x_1^5 + \cdots + x_5^5 )$. The enumerative theory for the LG phase was developed by Fan-Jarvis-Ruan~\cite{fjrw}, Chang-Kiem-Li~\cite{algvcquantumsingularity}, and Polishchuk-Vaintrob~\cite{mfcohftFJRW}.

A geometric construction for varying the K\"ahler parameters and connecting different phases of a GLSM is the theory of Mixed Spin $P$-fields, developed by Chang-Li-Li-Liu~\cite{mspfermat, msp2} based on the master space construction of Thaddeus~\cite{gitflips}. There are two key facts about MSP theory:
\begin{enumerate}
    \item It has a virtual localization formula which decomposes the theory into three parts: the geometric phase, Hodge integrals, and the LG phase;
    \item The total MSP theory behaves like the Gromov-Witten theory of a Fano variety, at least at genus $0$.
\end{enumerate}
Using the second fact, a key insight of Guo was to duplicate the field which corresponds to deformation of the K\"ahler parameter. This idea was realized in~\cite{nmsp} as $N$-Mixed Spin $P$-field theory\footnote{In this paper, we will refer to the theory with an arbitrary $N$ as simply MSP theory and treat it as depending on a parameter $N$. This is justified by the fact that we will consider the asymptotics of the theory at large $N$.} and was used by Chang-Guo-Li~\cite{nmsp2, nmsp3} to prove the Yamaguchi-Yau conjectures and the BCOV conjecture for the quintic threefold. 

In this paper, we will study the structure of the Gromov-Witten potentials $F_g$ for all $g \geq 0$ for the targets $Z_6 \subset \P(1,1,1,1,2)$, $Z_8 \subset \P(1,1,1,1,4)$, and $Z_{10} \subset \P(1,1,1,2,5)$ using MSP theory. Below, we will outline the results and methods used to prove them.

\subsection{Outline and statement of results}%
\label{sub:Outline and statement of results}

Let $\ba = (1,1,1,1,2)$, $(1,1,1,1,4)$, or $(1,1,1,2,5)$. Let $k = \sum_{i=1}^5 a_i$ and set
\[ p_k \coloneqq \int_{\P(\ba)} H^4 = \frac{k}{\prod_{i=1}^5 a_i}. \]

Define the $I$-function of $Z$ by
\begin{align*}
    I(q,z) &\coloneqq z\sum_{d \geq 0} q^d \frac{\prod_{m=1}^{kd}(kH+mz)}{\prod_{i=1}^5 \prod_{m=1}^{a_i d} (a_i H + mz)} \\
    &\eqqcolon I_0(q)z + I_1(q) H+ I_2(q) \frac{H^2}{z} + I_3(q) \frac{H^3}{z^2}.
\end{align*}
Note that the radius of convergence in $q$ of the $I$-function is $\frac{1}{r}$, where
\[ r \coloneqq \frac{k^k}{\prod_{i=1}^5 a_i^{a_i}}. \]

\begin{rmk}
    Our $I$-function is different from the usual $I$-function, which has a prefactor of $e^{\frac{H}{z}\log q}$. Therefore, applying the differential operator $zD \coloneqq zq \odv{}{q}$ to the usual $I$-function corresponds to applying the differential operator $(H+zD)$ to our $I$-function. 

    From the perspective of quasimap theory, our choice of convention for the $I$-function corresponds to the result obtained by localization on the stacky loop space~\cite[Proposition 4.9]{orbqmap} and its non-negative part is exactly the quasimap wall-crossing formula of~\cite{qmapwc}.
\end{rmk}

Let $D \coloneqq q\odv{}{q}$ and set 
\[I_{11} = 1+D\ab(\frac{I_1}{I_0}).\] 
Now define the infinite number of generators to be
\begin{align*}
    A_m \coloneqq \frac{D^m I_{11}}{I_{11}}, \qquad 
    B_m \coloneqq \frac{D^m I_{0}}{I_{0}}, \qquad \text{and} \qquad
    Y \coloneqq \frac{1}{1-rq}.
\end{align*}
Set $A \coloneqq A_1$ and $B \coloneqq B_1$. Define
\[ \msc{R} \coloneqq \Q[A,B, B_2, B_3, Y]. \]
\begin{lem}[\cite{yy04}\footnote{The proof of Yamaguchi-Yau relies on non-holomorphic completions of the generators and is written in the language of special geometry.}]
    The ring $\msc{R}$ contains all $A_m$ and $B_m$ for $m \geq 1$ and is closed under the operator $D$.
\end{lem}

Let $N_{g,d} = \ab<\ >_{g,0,d}^Z$ be the genus $g$ degree $d$ Gromov-Witten invariant of $Z$. We consider the generating function
\[ F_g(Z) = \begin{cases}
    a_{0,k} (\log Q)^3 + \sum_{d \geq 1} N_{0,d} Q^d & g=0 \\
    a_{1,k} \log Q + \sum_{d \geq 1} N_{1,d} Q^d & g=1 \\
    \sum_{d \geq 0} N_{g,d} Q^d & g\geq 2 ,
\end{cases}
\]
where $a_{0,k} = \int_Z H^3$ and $a_{1,k} = \int_Z c_2(Z) \cup H$ are given by~\Cref{tab:label}.
\begin{table}[htpb]
    \centering
    \caption{Values of $a_{0,k}$ and $a_{1,k}$}
    \label{tab:label}
    \begin{tabular}{ccc}
        \toprule
        $k$ & $a_{0,k}$ & $a_{1,k}$ \\
        \midrule 
        $k=6$ & $\frac{1}{2}$ & $-\frac{7}{4}$ \\
        \\[-0.8em]
        $k=8$ & $\frac{1}{3}$ & $-\frac{11}{6}$ \\
        \\[-0.8em]
        $k=10$ & $\frac{1}{6}$ & $-\frac{17}{12}$ \\
        \bottomrule
    \end{tabular}
\end{table}
Now introduce the normalized Gromov-Witten potential
\[ P_{g,n} \coloneqq \frac{(p_k Y)^{g-1} I_{11}^n}{I_0^{2g-2}} \ab(Q \odv{}{Q})^n F_g(Q) \Big\vert_{Q = qe^{\frac{I_1}{I_0}}}. \]
The main result of this paper is the finite generation conjecture of Yamaguchi-Yau:

\begin{thm}[\Cref{thm:polynomiality}]\label{thm:intropoly}
    For any genus $g$ and number of marked points $n$, we have
    \[ P_{g,n} \in \msc{R} = \Q[A,B,B_2, B_3, Y]. \]
\end{thm}

We will prove this theorem by studying the structure of MSP invariants. In this introduction, we give only an informal introduction to MSP theory, see~\Cref{sec:Setup and genus zero structure} for full details.

The moduli space $\msc{W}_{g,n,\bd}$ of MSP fields with numerical data $(g,n,\bd = (d_0, d_{\infty}))$ is defined in~\Cref{sub:invariants and state space}. The stack $\msc{W}_{g,n,\bd}$ has an action of $\T = (\C^{\times})^N$ and a $\T$-equivariant virtual cycle supported on a proper closed $\T$-substack $\msc{W}_{g,n,\bd}^-$. There are also $\T$-equivariant evaluation morphisms
\[ \ev_i \colon \msc{W}_{g,n,\bd} \to \P(\ba, 1^N). \]

If we consider the locus $\mf{C} = Z\ab(\sum_{i=1}^5 x_i^{\frac{k}{a_i}})$, then $\mf{C}$ has an action of $\T$. We will set the state space to be
\[ \msc{H} \coloneqq H^*(\mf{C}^{\T}, \Q(\zeta_N)(t)). \]
Because $\ev_i(\msc{W}_{g,n,\bd}^-) \subset \mf{C}$, we may define MSP correlators
\[ \ab<\tau_1(\psi_1), \ldots, \tau_n(\psi_n)>_{g,n,d_{\infty}}^M \]
using virtual localization. 

\begin{conv}\label{conv:eqvar}
    Let $\T = (\C^{\times})^N$. Then let the standard presentation of the $\T$-equivariant cohomology of a point be $H_{\T}(\mr{pt}, \Q) = \Q[t_1, \ldots, t_N]$. After equivariant integration, we will make the substitution $t_{\alpha} = - \zeta_N^{\alpha} t$, where $t$ is a formal variable of degree $1$. Finally, we will often make the change of variables $q' = \frac{-q}{t^N}$.
\end{conv}

Virtual localization in MSP theory splits the theory into three \textit{levels}, which are defined as follows:
\begin{itemize}
    \item Level $0$ yields the Gromov-Witten theory of $Z$;
    \item Level $1$ yields Hodge integrals;
    \item Level $\infty$ yields the LG phase of $Z$, and in particular, FJRW invariants.
\end{itemize}

We then construct a theory contained entirely in levels $0$ and $1$, called the $[0,1]$ theory, which coincides with the full MSP theory when $g=0$. In particular, we define classes
\[ [\tau_1(\psi_1), \ldots, \tau_n(\psi_n)]_{g,n}^{[0,1]} \in H^*(\ol{\msc{M}}_{g,n}) \]
and $[0,1]$ correlators
\[ \ab<\tau_1(\psi_1), \ldots, \tau_n(\psi_n)>_{g,n}^{[0,1]} = \int_{\ol{\msc{M}}_{g,n}} [\tau_1(\psi_1), \ldots, \tau_n(\psi_n)]_{g,n}^{[0,1]}. \]
The full MSP theory can be obtained from the $[0,1]$ theory and FJRW-like invariants at level $\infty$. They both satisfy the following degree bound:
\begin{thm}[\Cref{thm:polynomiality01},~\Cref{lem:fullpolynomiality}]\label{thm:intropoly01}
    Let $p \coloneqq c_1(\mc{O}_{\P(\ba, 1^N)}(1))$.
    Both the $[0,1]$ correlator
    \[ \ab<p^{a_1} \bar{\psi}_1^{m_1}, \ldots, p^{a_n}\bar{\psi}_n^{m_n}>_{g,n}^{[0,1]}. \]
    and the full MSP correlator
    \[ \ab<p^{a_1} \bar{\psi}_1^{m_1}, \ldots, p^{a_n}\bar{\psi}_n^{m_n}>_{g,n}^{M}. \]
    are both polynomials in the K\"ahler parameter $q$ of degree at most
    \[ g-1 + \frac{3g-3 + \sum_i a_i}{N}. \]
\end{thm}

We then relate it to the twisted Gromov-Witten theory of $\mc{C} \coloneqq \mf{C}^{\T} = Z \cup \bigcup_{\alpha=1}^N \pt_{\alpha}$ using the formalism of $R$-matrix actions on cohomological field theories:

\begin{thm}[\Cref{fthm:cohft},~\Cref{thm:01cohft}]\label{thm:cohftintro}
    The MSP $[0,1]$ theory for large $N$ defines a CohFT 
    \[\Omega^{[0,1]}(-) \coloneqq [-]^{[0,1]} \] which is related to the local theory
    \[ \Omega^{\loc}(-) \coloneqq [-]^{\loc} \]
    defined in~\Cref{eqn:localcohft}
    by the action of the $R$-matrix defined in~\Cref{eqn:rmatrix} by
    \[ \Omega^{[0,1]} = R.\Omega^{\loc}. \]
\end{thm}

This means that the $[0,1]$ theory can be obtained as an explicit sum of contributions over stable graphs. After proving polynomiality of the entries of the $R$-matrix in~\Cref{sec:Computation of the R-matrix}, we combine them with the stable graph localization formula to obtain~\Cref{thm:intropoly}. In the process of the proof, we directly compute $P_{1,1}$ and obtain the following formulae~\Cref{eqn:p11}:

\begin{thm}
    We have the formula    
    \[
    P_{1,1} = -\frac{1}{2} A + \ab(\frac{\chi(Z)}{24} - 2)B - \frac{1}{12}X + a_{1,k}, 
    \]
    where $\chi(Z_6) = -204$, $\chi(Z_8) = -296$, $\chi(Z_{10}) = -288$, and $a_{1,k} = \int_{Z_k} c_2(Z_k) \cup H$ is given in~\Cref{tab:label}.
\end{thm}

This paper is organized as follows:
\begin{itemize}
    \item In~\Cref{sec:Setup and genus zero structure}, we give a precise introduction to MSP theory, introduce the $[0,1]$ theory, and use the Givental formalism to study the genus zero MSP theory.
    \item In~\Cref{sec:MSP theory as a CohFT}, we package the $[0,1]$ theory using Givental's quantization formalism and prove~\Cref{thm:cohftintro}.
    \item In~\Cref{sec:Polynomiality of the 01 theory}, we explain how to obtain decompose the full MSP theory into $[0,1]$ and $(1\infty]$ parts and prove~\Cref{thm:intropoly01}.
    \item In~\Cref{sec:Computation of the R-matrix}, we prove the polynomiality of the entries of the $R$-matrix;
    \item In~\Cref{sec:Structure of the GW potential}, we use the results from the previous sections to prove~\Cref{thm:intropoly}.
\end{itemize}
In subsequent work~\cite{bcovme}, we will build on these results to prove the BCOV conjecture for the targets $Z_6 \subset \P(1,1,1,1,2)$, $Z_8 \subset \P(1,1,1,1,4)$, and $Z_{10} \subset \P(1,1,1,2,5)$.

\subsection*{Acknowledgements}%
\label{sub:Acknowledgements}

The author would like to thank Chiu-Chu Melissa Liu for proposing the use of MSP fields to study the Gromov-Witten theory of Calabi-Yau hypersurfaces in toric orbifolds and for explaining virtual localization, Dimitri Zvonkine for his lectures about cohomological field theories and Teleman's theorem at the Simons Center for Geometry and Physics in August 2023, and Konstantin Aleshkin and Shuai Guo for helpful discussions. The author would finally like to thank Jun Li and Shuai Guo for their hospitality during the author's visits to the Shanghai Center for Mathematical Sciences in May 2024 and to Peking University in July 2024, respectively, where part of this work was completed.

\section{Setup and genus zero structure}%
\label{sec:Setup and genus zero structure}

This section sets up MSP theory for our targets and computes the quantum connection and some properties of its fundamental solution. For more details, see~\cite{foundations}.

\subsection{MSP invariants and state space}%
\label{sub:invariants and state space}

Let $\msc{W}_{g,n,\mbf{d}}$ be the moduli of MSP fields, where all marked points have type $(1,\rho)$. Recall that a closed point $\xi \in \msc{W}_{g,n,\mbf{d}}$ is a tuple
\[ \xi = (\msc{C}, \Sigma^{\msc{C}}, \msc{L}, \msc{N}, \varphi , \rho , \mu , \nu ) \]
where $(\msc{C}, \Sigma^{\msc{C}})$ is a prestable twisted curve, $\msc{L}, \msc{N}$ are line bundles on $\msc{C}$ such that $(\msc{L}, \msc{N})$ is representable (as a morphism to $\mbf{B}\G_m^2$), and the fields are
\[ \varphi \in \bigoplus_{i=1}^5 H^0(\msc{L}^{a_i}),\ \rho \in H^0(\msc{L}^{-k} \otimes \omega_{\msc{C}}^{\log}),\ \mu \in H^0(\msc{L} \otimes \msc{N})^{\oplus N},\ \text{and}\ \nu \in H^0(\msc{N}). \]
We require that $\Aut(\xi)$ is finite and that the fields satisfy the properties that $(\varphi, \mu)$, $(\rho, \nu)$, and $(\mu, \nu)$ are nowhere vanishing. In addition, we require that $\rho$ vanishes at the marked points. The discrete data we will consider are the genus $g$, the number of markings $n$, and the degrees $d_0 = \deg (\msc{L} \otimes \msc{N})$ and $d_{\infty} = \deg \msc{N}$. 

The moduli $\msc{W}_{g,n,\mbf{d}}$ is a separated Deligne-Mumford stack which has an action of $\T = (\C^{\times})^N$ by scaling $\mu$. It has a perfect obstruction theory equipped with a cosection that has proper degeneracy locus $\msc{W}_{g,n,\mbf{d}}^-$, so we obtained a $\T$-equivariant virtual cycle $[\msc{W}_{g,n,\mbf{d}}]^{\vir}$.

Note that there is an evaluation map $\ev_j \colon \msc{W}_{g,n,\mbf{d}} \to \P(\mbf{a}, 1^n)$ given by
\[ \xi \mapsto \ab[\varphi(\Sigma_j), \frac{\mu(\Sigma_j)}{\nu(\Sigma_j)}] \]
which exists because $\nu$ is nowhere vanishing at the markings. If we allow $\T$ to scale the last $N$ coordinates of $\P(\mbf{a}, 1^N)$, then this morphism is clearly $\T$-equivariant.

Note that virtual localization gives
\[ [\msc{W}_{g,n,\bd}]^{\vir} = \sum_{\Theta \in G_{g,n,\bd}^{\mr{reg}}} \frac{[F_{\Theta}]^{\vir}}{e(N_{\Theta}^{\vir})}, \]
where the sum is over the set $G_{g,n,\bd}^{\mr{reg}}$ of \textbf{regular} graphs. Recall that these have vertices $V$, edges $E$, and legs $L$ together with decorations. The decorations are as follows:
\begin{enumerate}
    \item Every vertex is assigned a \textit{level}, which is either $0$, $1$, or $\infty$; Edges are partitioned into $E = E_{01} \cup E_{11} \cup E_{1\infty} \cup E_{\infty\infty}$ based on the levels of the two incident vertices. Because $\Theta$ is regular, $E_{0\infty}$ is empty.
    \item Every vertex in $V_1 \cup V_{\infty}$ is assigned an \textit{hour} in $\{1, \ldots, N\}$.
    \item Every vertex is assigned a genus and degree.
\end{enumerate}

Recall that $\msc{W}_{g,n,\bd}^- \subset \msc{W}_{g,n,\bd}$ is the locus where either $\varphi = 0$ or $\rho = \sum_{i=1}^5 \varphi_i^{\frac{k}{a_i}} = 0$. In particular, because $\rho = 0$, we see that
\[ \ev_j(\msc{W}_{g,n,\bd}^-) \subset \ab(\sum_i x_i^{\frac{k}{a_i}} = 0) \subset \P(\ba, 1^N). \]
The fixed locus $\P(\ba, 1^N)$ is a union of $\P(\ba)$ at level $0$ and $N$ isolated points $N\pt = \bigcup_{\alpha=1}^N \ab\{\pt_{\alpha}\}$ at level $1$. Because markings in $F_{\Theta}$ are attached to vertices of $\Theta$, we see that
\[ \ev_j(F_{\Theta} \cap \W_{g,n,\bd}^-) \subset Z \cup N\pt \eqqcolon \mc{C}. \]

We will now introduce the state space. Let $\F \coloneqq \Q(\zeta_N)(t)$. By~\Cref{conv:eqvar}, all MSP invariants will lie in $\F$ and their generating series will lie in $\A \coloneqq F\llbracket q \rrbracket$.

\begin{defn}
    Define the \textit{state space} to be
    \[ \msc{H} = H^*(\mc{C}, \F). \]
    Denote its even part by $\msc{H}^{\ev}$.
\end{defn}

Note that there is a surjective ring homomorphism
\[ H_{\T}^*(\P(\ba, 1^N), \Q)|_{t_{\alpha} \mapsto -\zeta_N^{\alpha}t, \loc} = \Q(t)[p]\Bigg/\ab(p^5 \prod_{\alpha=1}^N (p+t_{\alpha})) \to \msc{H}^{\ev} \]
such that $p|_{Z} = H$ and $p|_{\pt_{\alpha}} = -t_{\alpha}$, where $p = c_1(\msc{O}_{\P(\ba, 1^N)}(1))$. Its kernel is spanned by $p^4 \cdot \prod_{\alpha=1}^N (p+t_{\alpha})$. By Lagrange interpolation, we see that
\begin{align*}
    \mbf{1}_{\alpha} &= \frac{p^4}{t_{\alpha}^4} \prod_{\beta \neq \alpha} \frac{p+t_{\beta}}{t_{\beta} - t_{\alpha}} \\
    H^j &= \frac{p^j}{t^N} (t^N - p^N)
\end{align*}
for all $\alpha = 1, \ldots, N$ and $j = 0, \ldots, 3$.

We may finally define the \textit{MSP invariants}
\[ \ab<\bigotimes_{j=1}^n \tau_j(\psi_j)>_{g,n,d_{\infty}}^M \coloneqq (-1)^{1-g}\sum_{d \geq 0} q^d \sum_{\Theta \in G_{g,n,(d,d_{\infty})}^{\mr{reg}}} \int_{[F_{\Theta}]^{\vir}}\frac{\prod_{j=1}^N \ev_j^*\tau_j(\psi_j)}{e(N_{\Theta}^{\vir})}. \]
By the virtual localization formula, this is equal to
\[ \sum_{d \geq 0}q^d \int_{[\msc{W}_{g,n,(d,d_{\infty})}]^{\vir}} \prod_{j=1}^N \ev_j^*\tau_j(\psi_j) \]
whenever we consider ambient insertions. When $d_{\infty} = 0$, we abbreviate
\[ \ab<->_{g,n}^M \coloneqq \ab<->_{g,n,0}^M. \]
The following lemma will help us understand the MSP invariants in genus $0$.

\begin{lem}\label{lem:genus0moduli}
    When $g=0$ and $d_{\infty} = 0$, there is an isomorphism
    \[ \msc{W}_{0,n,(d, 0)} \cong \ol{\msc{M}}_{0,n}(\P(\mbf{a}, 1^N), d). \]
\end{lem}

\begin{proof}
    First, we show that $\nu$ is nowhere vanishing. Let $\msc{E} \subseteq \msc{C}$ be a connected component of $(\nu = 0)$. Then we note that $\mu|_{\msc{E}}$ and $\rho|_{\msc{E}}$ are nowhere vanishing, so $\msc{E}$ contains no markings, $\msc{L}^{-k} \otimes \omega_{\msc{C}}^{\log}|_{\msc{C}} \cong \msc{O}_{\msc{E}}$, and $\deg \msc{L} \otimes \msc{N}|_{\msc{E}} \geq 0$. If $m$ is the number of nodes (connecting $\msc{E}$ to irreducible components of $\msc{C}$ not contained in $\msc{E}$) attached to $\msc{E}$, then we see that $\deg \msc{L}|_{\msc{E}} = \frac{m-2}{k}$. However, each irreducible component attached to $\msc{E}$ has $\deg \msc{N} \geq \frac{1}{k}$, so $\msc{E}$ contributes at least $\frac{2}{k} > 0$ to $\deg \msc{N} = d_{\infty}$. All components of $\msc{C}$ not intersecting $(\nu = 0)$ have $d_{\infty} = 0$, so we see that if $(\nu = 0) \neq \emptyset$, then $d_{\infty} > 0$.

    We now see that $\nu$ is nowhere vanishing, so $\msc{N} \cong \msc{O}_{\msc{C}}$. Because $g = 0$ and $\deg \msc{L} \geq 0$, we know that $\deg(\msc{L}^{-k} \otimes \omega_{\msc{C}}^{\log}) \leq n-2$. Because $\rho$ must vanish at the $n$ marked points, it is identically zero.
\end{proof}

Note that~\Cref{lem:genus0moduli} tells us that the genus zero MSP theory equals the (untwisted part of the) genus zero $T$-equivariant $\msc{O}(k)$-twisted Gromov-Witten theory of $\P(\ba, 1^N)$, or equivalently, the ambient part of the genus-zero Gromov-Witten theory of a degree $k$ hypersurface in $\P(\ba, 1^N)$. This implies that the invariants satisfy the string and dilaton equations as well as the topological recursion relations.

Assuming $N > 3$, we define the pairing
\begin{align*}
    (x,y)^M &\coloneqq [q^0]\ab<1,x,y>_{0,3}^M \\
    &= \int_{Z} \frac{xy}{\prod_{\alpha=1}^N (H+t_{\alpha})} + \sum_{\alpha} \frac{p_k xy}{t_{\alpha}^4 \prod_{\beta \neq \alpha} (t_{\beta} - t_{\alpha})} \bigg\vert_{\pt_{\alpha}} \\
    &= \int_{Z} \frac{xy}{-t^N} + \sum_{\alpha} \frac{-p_k}{N t_{\alpha}^3 t^N} xy \bigg\vert_{\pt_{\alpha}}\\
    &\eqqcolon (x|_{Z},y|_{Z})^{Z, \tw} + \sum_{\alpha} (x|_{\pt_{\alpha}}, y|_{\pt_{\alpha}})^{\pt_{\alpha}, \tw}.
\end{align*}
In particular, using~\Cref{lem:genus0moduli}, in the basis $\{ \phi_i \coloneqq p^i\}$, the pairing is given by
\begin{align*}
    \eta_{ij} &\coloneqq (p^i, p^j)^M \\
    &= \int_{\P(\ba, 1^N)} kp \cup p^i \cup p^j \\
    &= \begin{cases}
        p_k t^N & i+j = 2N+3 \\
        p_k & i+j = N+3 \\
        0 & \text{otherwise}.
    \end{cases}
\end{align*}
In this basis, the Poincar\'e duals are given by
\[ \phi_j = \frac{p^{3-j}}{p_k}(p^N-t^N) \]
for $j=0, \ldots, 3$ and
\[ \phi^j = \frac{p^{N+3-j}}{p_k} \]
for $j=4, \ldots, N+3$. In the $\{H^j\} \cup \{\mbf{1}_{\alpha}\}$ basis, the Poincar\'e duals are given by
\[ \ab\{ -\frac{t^N}{p_k} H^{3-j}\} \cup \ab\{\mbf{1}^{\alpha} \coloneqq \frac{N t_{\alpha}^3 t^n}{p_k} \mbf{1}_{\alpha} \}. \]

\subsection{The $[0,1]$ and local theories}
\label{sub:The 01 and local theories}

We will now define the MSP $[0,1]$ theory as well as the local theory of $\mc{C}$. The local theory corresponds to the local contributions of vertices in the virtual localization formula for the MSP moduli space. The two theories will be related by an $R$-matrix in a way akin to Teleman's classification of semisimple cohomological field theories~\cite{2dsscohft}.
\begin{defn}
    We define the \textit{twisted classes} for $Z$ by the formula 
    \[ [\ol{\msc{M}}_{g,n}(Z, d)]^{\tw} \coloneqq \prod_{\alpha} e_{\T}(R \pi_* f^* \msc{O}(1) \cdot t_{\alpha})^{-1} \cap [\ol{\msc{M}}_{g,n}(Z, d)]^{\vir} \]
    and the twisted classes for $\pt_{\alpha}$ by
    \begin{align*}
        &[\ol{\msc{M}}_{g,n}]^{\alpha, \tw} \\ 
        \coloneqq{} &(-1)^{1-g} \frac{p_k t_{\alpha} \cdot \prod_{i=1}^5 e_{\T}(\E_{g,n}^{\vee}\cdot (-a_i t_{\alpha})) \cdot \prod_{\beta \neq \alpha} e_{\T}(\E_{g,n}^{\vee} \cdot (t_{\beta} - t_{\alpha}))}{(-t_{\alpha})^5 \cdot e_{\T}(\E_{g,n} \cdot kt_{\alpha}) \cdot \prod_{\beta \neq \alpha} (t_{\beta} - t_{\alpha})} \cap [\ol{\msc{M}}_{g,n}].
    \end{align*}
\end{defn}

In computations, we will only see the top degree parts of these twisted classes, which we will denote by $[\ol{\msc{M}}_{g,n}(Z, d)]^{\top}$ and $[\ol{\msc{M}}_{g,n}]^{\top}$.
\begin{lem}\label{lem:topclasses}
    Assuming $N > 3$, we have
    \begin{align*}
        [\ol{\msc{M}}_{g,n}(Z, d)]^{\top} &= [\ol{\msc{M}}_{g,n}(Z, d)]^{\tw} \\
        &= (-t^N)^{-(d+1-g)} [\ol{\msc{M}}_{g,n}(Z, d)]^{\vir}
    \end{align*}
    for the target $Z$. For the isolated points $\pt_{\alpha}$, we have
    \[ [\ol{\msc{M}}_{g,n}]^{\top} = \ab(\frac{1}{p_k} N(-t_{\alpha})^{N+3})^{g-1} [\ol{\msc{M}}_{g,n}]. \]
\end{lem}
The proof of the above lemma uses the same argument used to prove~\cite[Lemma 5.9]{msp2}.

We will now define local invariants and classes.
\begin{defns}
    For $\tau_i(\psi) = \sum_j \tau_{i,j} \psi^j \in \msc{H}[\psi]$, we define
    \begin{align*}
        \ab<\bigotimes_{i=1}^n \tau_i(\psi_i)>_{g,n}^{Z, \tw} &\coloneqq \sum_{d \geq 0} q^d \int_{[\ol{\msc{M}}_{g,n}(Z,d)]^{\tw}} \prod_i \ev_i^*(\tau_i(\psi_i)|_{Z}) \\
        \ab<\bigotimes_{i=1}^n \tau_i(\psi_i)>_{g,n}^{\pt_{\alpha}, \tw} &\coloneqq \int_{[\ol{\msc{M}}_{g,n}]^{\alpha, \tw}} \prod_i (\tau_i(\psi_i)|_{\pt_{\alpha}}).
    \end{align*}
    Deleting the ``tw'' defines $\ab<\bigotimes_i \tau_i(\psi_i)>_{g,n}^{Z}$ and $\ab<\bigotimes_i \tau_i(\psi_i)>_{g,n}^{\pt_{\alpha}}$. Now define
    \[ \ab<->_{g,n}^{\loc} \coloneqq \ab<->_{g,n}^{Z, \tw} + \sum_{\alpha} \ab<->_{g,n}^{\pt_{\alpha}, \tw}. \]

    When $2g-2+n > 0$, consider the diagram
    \begin{equation*}
    \begin{tikzcd}
        \ol{ \msc{M} }_{g,n}(Z, d) \ar{r}{\ev_i} \ar{d}{\st^{Z}} & Z \\
        \ol{\msc{M}}_{g,n}.
    \end{tikzcd}
    \end{equation*}
    Now define
    \begin{align*}
        \ab[\bigotimes_{i=1}^n \tau_i(\psi_i)]^{Z, \tw}_{g,n} &\coloneqq \sum_{d \geq 0} q^d \st^{Z}_* \ab(\prod_i \ev_i^*(\tau_i(\psi_i)|_{Z}) \cap [\ol{\msc{M}}_{g,n}(Z, d)]^{\tw}) \\
        \ab[\bigotimes_{i=1}^n \tau_i(\psi_i)]^{\pt_{\alpha} \tw}_{g,n} &\coloneqq \prod_i (\tau_i(\psi_i)|_{\pt_{\alpha}}) \cap [\ol{\msc{M}}_{g,n}]^{\alpha, \tw}.
    \end{align*}
    Deleting ``tw'' defines $[-]_{g,n}^{Z}$ and $[-]_{g,n}^{\pt_{\alpha}}$ and replacing ``tw'' with ``top'' defines $[-]_{g,n}^{Z, \top}$ and $[-]_{g,n}^{\pt_{\alpha}, \top}$. We finally define
    \begin{equation}\label{eqn:localcohft}
    [-]_{g,n}^{\loc} \coloneqq [-]_{g,n}^{Z, \top} + \sum_{\alpha} [-]_{g,n}^{\pt_{\alpha}, \top}. 
    \end{equation}
\end{defns}

We will now define the MSP $[0,1]$ theory.
\begin{defn}
    Let $G_{g,n,d}^{[0,1]}$ denote the set of graphs in $G_{g,n,(d,0)}^{\mr{reg}}$ where all vertices have level $0$ or $1$. Define the $[0,1]$-virtual cycle by
    \[ [W_{g,n,(d,0)}]^{[0,1]} = \sum_{\Theta \in G_{g,n,d}^{[0,1]}} \frac{[F_{\Theta}]^{\vir}}{e(N_{\Theta}^{\vir})}. \]
\end{defn}

\begin{defn}
    Let $g,n$ be such that $2g-2+n > 0$. For any $\tau_i(z) \in \msc{H}\ps{z}$, define
    \[\ab[\bigotimes_{i=1}^n \tau_i(\psi_i)]_{g,n}^{[0,1]} \coloneqq \sum_{d \geq 0} (-1)^{1-g} q^d \st^{M}_* \ab(\prod_{i=1}^n \ev_i^* \tau_i(\psi_i) \cap [\msc{W}_{g,n,(d,0)}]^{[0,1]}), \]
    where $\st^M \colon \msc{W}_{g,n,(d,0)} \to \ol{\msc{M}}_{g,n}$ stabilizes the coarse moduli space of $(\msc{C}, \Sigma^{\msc{C}})$.

    We also define the MSP $[0,1]$ invariants by
    \[ \ab<\bigotimes_{i=1}^n \tau_i(\psi_i)>_{g,n}^{[0,1]} \coloneqq \int_{\ol{\msc{M}}_{g,n}} \ab[\bigotimes_{i=1}^n \tau_i(\psi_i)]_{g,n}^{[0,1]}. \]
\end{defn}

\begin{rmk}
    In genus $0$, the $[0,1]$ invariants are the equal to the full MSP invariants by~\Cref{lem:genus0moduli}.
\end{rmk}

\subsection{Genus zero MSP theory}%
\label{sub:genus zero msp}

Let ``$\bullet$'' be either ``$M$'' or ``loc'' and choose a basis $\{e_{\alpha}\}$ of $\msc{H}$. Denote its dual basis under $(-,-)^M$ by $\{e^{\alpha}\}$. We define the $J$-function and $S$-matrix by the formulae
\begin{align}
    J^{\bullet}(\mbf{t}, z) &\coloneqq z + \mbf{t}(-z) + \sum_{\alpha, n} \frac{1}{n!} e^{\alpha} \ab<\frac{e_{\alpha}}{z-\psi}, \mbf{t}(\psi)^n>_{0,n+1}^{\bullet}, \label{eqn:jfunction} \\
    S_{\tau}^{\bullet}(z)x &\coloneqq x + \sum_{\alpha, n} \frac{1}{n!} e^{\alpha} \ab<\frac{x}{z-\psi},e_{\alpha}, \tau^n>_{0,n+2}^{\bullet}  \label{eqn:smatrix}
\end{align}
for $\mbf{t} \in \A \otimes \msc{H}\ps{z}$ and $\tau \in \A \otimes \msc{H}$.

\begin{conv}
    We will expand $J^M$ and $S^M_{\tau}$ near $z=0$ and expand $J^{\loc}$ and $S^{\loc}_{\tau}$ near $z=\infty$.
\end{conv}

\begin{rmks}
    We make some observations about the $J$-functions and $S$-matrices.
    \begin{enumerate}
        \item By degree reasons, $S_{\tau}^M$ preserves $\msc{H}^{\ev}$ and $H^3(Z)$ whenever $\tau \in \msc{H}^{\ev}$. We will now consider $x \in H^3(Z)$. The same is true for $S^{\loc}_{\tau}$. 
        \item For each $d$, the coefficient of $q^d$ in $J^{M}(\mbf{t}, z)$ is a rational function in $z$ whenever $\mbf{t} \in \A \otimes \msc{H}^{\geq 2}\ps{z}$. The reason is that when we compute $J^M$ by localization, there are finitely many localization graphs. For each graph, either the insertion $\frac{e_{\alpha}}{z-\psi}$ is attached to a stable or unstable vertex. If the vertex is stable, then for degree reasons only a finite number of powers of $\psi$ will give nonzero invariants, and if the vertex is unstable, we obtain something like $\frac{1}{z-w}$, where $w$ is some equivariant weight. The same holds for $S_{\tau}^{M}(z) x$ for any $x \in \msc{H}$ and $ \tau \in \A \otimes \msc{H}^{\geq 2}$.
        \item With the same assumptions $\mbf{t} \in \A \otimes \msc{H}^{\geq 2} \ps{z}$ and $\tau \in \A \otimes \msc{H}^{\geq 2}$, the same conclusions hold for $J^{\loc}(\mbf{t}, z)$ and $S_{\tau}^{\loc}(z)x$ for any $x \in \msc{H}$. This is because any element of $\msc{H}_{\geq 2}$ is supported on $Z$ and $J^Z$ and $S^Z_{\tau}$ can be computed by localization on $\P(\ba)$.
        \item When $x \in H^3(Z)$, we compute $S_{\tau=0}^M(z)x = x$ by localization. First, we note that both the $x$ and $e_{\alpha}$ insertions must be at level $0$, and then for dimension reasons the contributions of the vertices containing the markings must be $0$. In addition, for $\tau \in \msc{H}^{\ev}$ such that $\tau|_Z = \sum_{i=0}^3 \tau_i H^i$ and $x \in H^3(Z)$, the string equation and a virtual dimension computation yields $S_{\tau}^{\loc}(z)x = e^{\frac{\tau_0}{z}} x$.
    \end{enumerate}
\end{rmks}

The same formulae for $F = Z$ or $\pt_{\alpha}$ and $\bullet$ being ``$F$'' or ``$F, \tw$'' define the $J$-functions and $S$-matrices $J^{\bullet}(\mbf{t}, z) \in \A \otimes \msc{H}_F\ps{z,z^{-1}}$ and $S_{\tau}^{\bullet}(z) \in \A \otimes \End \msc{H}_F\ps{z^{-1}}$, where $\msc{H}_F = H^*(F, \F)$. Here, we take $\tau \in \A \otimes \msc{H}_F$ and $\mbf{t} \in \A \otimes \msc{H}_F\ps{z}$. The pairings for $\bullet$ being ``$F$'' are the ordinary Poincar\'e pairings while the the pairings for $\bullet$ being ``$F,\tw$'' are the twisted Poincar\'e pairings $(-,-)^{F, \tw}$ defined previously.

By direct computation, we obtain the formulae
\begin{align*}
    S_{\tau}^{Z, \tw}(q,z) &= S_{\tau}^{Z}(q', z), \\
    S_{\tau}^{\pt_{\alpha}, \tw} &= S_{\tau}^{\pt_{\alpha}}(z) = e^{\frac{\tau}{z}}.
\end{align*}

By~\Cref{lem:genus0moduli}, we may apply the mirror theorem of Coates-Corti-Iritani-Tseng for complete intersections in toric stacks.
\begin{thm}[{\cite[Corollary 28]{applmirrortoricstack}}]\label{thm:mirror}
    Consider the (very small) MSP $J$-function $J^M(z) \coloneqq J^M(0,z)$ and define the MSP $I$-function
    \begin{equation}\label{eqn:mspifn}
        I^M(q, z) = z + \sum_{d=1}^{\infty} z q^d \frac{\prod_{m=1}^{kd}(kp+mz)}{\prod_{i=1}^5 \prod_{m=1}^{a_i d}(a_i p + mz) \prod_{m=1}^d ((p+mz)^N - t^N)}.
    \end{equation}
    Then whenever $N > 1$, we have
    \[ J^M(z) = I^M(q, z). \]
\end{thm}

We will now relate $J^M$ and $J^{\loc}$. The main tool for this is the Lagrangian cone, which was introduced in~\cite{qrr,symplfrob}. For any (linear) Frobenius manifold with underlying vector space $\msc{H}$ and potential $F_0$ (in Gromov-Witten theory, this is the genus zero descendent potential), we introduce the vector space $\msc{H}\llparenthesis z^{-1} \rrparenthesis$ with the symplectic form
\[ \Omega(f,g) = \on{Res}_{z=0}(f(-z), g(z)). \]
This polarizes the symplectic space into non-negative and negative powers of $z$. We then choose Darboux coordinates $(\mbf{p}, \mbf{q})$ on $\msc{H}(\llparenthesis z^{-1} \rrparenthesis)$ and define the \textit{Lagrangian cone} to be the graph of $\d F_0$, or in other words,
\[ \msc{L} \coloneqq \ab\{(\mbf{p}, \mbf{q}) \coloneqq \mbf{p}(-z) + q(z) \in \A \otimes \msc{H}\llparenthesis z^{-1} \rrparenthesis \mid \mbf{p} = \d_{\mbf{q}} F_0(\mbf{t}) \}. \]
Here, $\mbf{q}(z) = \mbf{t}(z) - z$; this is called the \textit{dilaton shift}.

\begin{thm}[{\cite[Theorem 1]{symplfrob}}]
    $\msc{L}$ is a Lagrangian cone with vertex at $\mbf{q} = 0$ such that its tangent spaces $L$ are tangent to $\msc{L}$ exactly at $zL$.
\end{thm}

Applying this theory to
\[ F_0^{\loc}(\mbf{t}) \coloneqq F_g^{Z, \tw}(\mbf{t}|_{Z}) + \sum_{\alpha} F_g^{\pt_{\alpha}, \tw}(\mbf{t}|_{\pt_{\alpha}}), \]
we obtain the Lagrangian cone $\msc{L}^{\loc}$ for the local theory.

\begin{rmk}
    An important property of the Lagrangian cone is that it can be recovered from the $J$-function.
\end{rmk}

\begin{lem}\label{lem:msponlocallagcone}
    For all $\mbf{t} \in \A \otimes \msc{H}^{\ev}\ps{z}$ and $\ep, \tau \in \A \otimes \msc{H}^{\ev}$, we have $J^M(\mbf{t}, -z) \in \msc{L}^{\loc}$ and $z S_{\ep}^M(z)^{-1} \tau \in z S_{\tau}^{\loc}(z)^{-1} \A \otimes \msc{H}\ps{z}$.
\end{lem}

\begin{proof}
    The proof is the same as~\cite[Lemma 1.8]{nmsp2}.
\end{proof}

\begin{notn}
    For a power series $f$ in $z, z^{-1}$, we will use $f(z^-)$ to denote expanding $f(z)$ around $z = \infty$ and $[f]_+$ to denote the part of $f$ with non-negative powers of $z$.
\end{notn}

In the proof of~\cite[Lemma 1.8]{nmsp2}, the non-negative part
\[ L^{\loc}(\mbf{t}, z) \coloneqq z \1 + [J^M(\mbf{t}, -z)]_+ \in \msc{H} \otimes \A\ps{z} \]
was introduced. This records the contributions of rational tails leading to unstable vertices, so we have
\[ J^M(\mbf{t}, -z) = J^{\loc}(L^{\loc}, -z). \]

\begin{cor}
    There exists a unique $R_{\ep}^{\loc}(z) \in \End \msc{H} \otimes \A\ps{z}$ such that
    \begin{equation}\label{eqn:rloc}
    S_{\ep}^M(z) = R_{\ep}^{\loc}(z) S_{\tau^{\loc}}^{\loc}(z), 
    \end{equation}
    where $\tau^{\loc} = \tau^{\loc}(\ep, q)$ is defined by the Dijkgraaf-Witten map\footnote{The terminology follows that used in~\cite{virasorotoricbundle,nmsp2}; the formula is originally due to~\cite{dijkgraafwitten}.}
    \[ \tau^{\loc}(\ep, q) \coloneqq \sum_{\alpha, n} \frac{1}{n!} e^{\alpha} \ab<e_{\alpha}, \mbf{1}, L^{\loc}(\ep, -\psi)^n>_{0,n+2}^{\loc} \in \msc{H} \otimes \A. \]
\end{cor}

We will use $\tau_{Z}(q') \coloneqq \tau^{\loc}(0,q)|_{Z}$ and $\tau_{\alpha}(q') \coloneqq \tau^{\loc}(0,q)|_{\pt_{\alpha}}$ for the components of the Dijkgraaf-Witten map. Later, we will compute explicit formulae
\begin{align*}
    \tau_{\alpha}(q) &= \zeta_N^{\alpha}t \int_{0}^q ( (1-rx)^{\frac{1}{N}}-1 ) \frac{\d x}{x}, \\
    \tau_{Z}(q) &= \frac{I_1(q)}{I_0(q)}H.
\end{align*}

\begin{conv}
    For the rest of this paper, we will abbreviate 
    \begin{align*}
        &S^M \coloneqq S_{\ep=0}^M, && S^{\loc} \coloneqq S^{\loc}_{\tau^{\loc}(0,q)}, && R^{\loc} \coloneqq R^{\loc}_{\ep=0}, && S^{Z, \tw} \coloneqq S^{Z,\tw}_{\tau_{Z}(q')}, \\ 
        &S^{Z} \coloneqq S^{Z}_{\tau_{Z}(q)}, && S^{\pt_{\alpha}, \tw} \coloneqq S_{\tau_{\alpha}(q')}^{\pt_{\alpha},\tw},  && \text{and } && S^{\pt_{\alpha}} \coloneqq S_{\tau_{\alpha}(q)}^{\pt_{\alpha}}. 
    \end{align*}
\end{conv}

We will now give an explicit formula for the MSP $S$-matrix on $\msc{H}^{\ev}$, again using~\Cref{lem:genus0moduli} to see that the Lagrangian cone $\msc{L}^M$ is the same as the Lagrangian cone of the $\msc{O}(k)$-twisted Gromov-Witten theory of $\P(\ba, 1^N)$ (ignoring the twisted sector).
\begin{lem}\label{lem:quantumconnection}
    Define $D_p \coloneqq p + z q \odv{}{q}$. Restricted to $\msc{H}^{\ev}$, $S^M(z)^*$ satisfies the quantum differential equation
    \[ D_p S^M(z)^* = S^M(z)^* \cdot A^M, \]
    where $A^M \in \End \msc{H}^{\ev}$ is given by the following matrix in the basis $\{\phi_i = p^i\}$ for all $N \geq 5$:
    \[ \begin{bmatrix}
        0 & & & & & & &  c_{1,k} q \\
        1 & 0 & & & & & & & c_{2,k} q \\
        & 1 & 0 & & & & & & & c_{3,k} q \\
        & & 1 & 0 & & & & & & & c_{2,k} q \\
        & & & 1 & 0 & & & & & & & c_{1,k} q + t^N \\
        & & & & 1 & 0 \\
        & & & & & 1 & 0 \\
        & & & & & & 1 & 0 \\
        & & & & & & & \cdots & \cdots \\
        & & & & & & & & 1 & 0 \\
        & & & & & & & & & 1 & 0 \\
        & & & & & & & & & & 1 & 0 \\
    \end{bmatrix}.
    \]
    Here, the coefficients in the top right corner are given by $c_{*,6} = (360,2772,5400)$, $c_{*,8} = (1680, 15808, 30560)$, and $c_{*,10} = (15120, 179520, 410720)$.
\end{lem}

\begin{proof}
    We compute $S^M(z)^* \phi_i$ recursively following the proof of~\cite[Lemma 1.12]{nmsp2}.
\end{proof}

Later, we will need the following specializations of the $S$-matrix.
\begin{defn}\label{defn:specializeds}
    Define
    \begin{align*}
        \msc{S}_a^{\alpha} &\coloneqq S^M(z) \mbf{1}^{\alpha}|_{z = \frac{k t_{\alpha}}{a}} \\
        \msc{S}_{a;i}^{\alpha} &\coloneqq (\msc{S}_a^{\alpha}, p^i)^M.
    \end{align*}
\end{defn}

\begin{lem}
    The following explicit formula for $\msc{S}_{a;0}^{\alpha}$ holds:
    \[ \msc{S}_{a;0}^{\alpha} = 1 + \sum_{d = 1}^{\ceil{a/k}-1} q^d \frac{(a-1)_{kd} \ab(\frac{a}{kt})^{Nd}}{\prod_{i=1}^5\ab(\frac{a_i a}{k}-1)_{a_i d} \prod_{m=1}^d \ab(\ab(-\frac{a}{k}+m)^N - \ab(\frac{a}{k})^N)} , \]
    where $(a)_m \coloneqq a(a-1) \cdots (a-m+1)$ is the falling Pochhammer symbol.
\end{lem}

\begin{proof}
    We compute
    \begin{align*}
        \frac{k t_{\alpha}}{a} \msc{S}_{a;0}^{\alpha} \mbf{1}_{\alpha} &= \frac{k t_{\alpha}}{a} (\msc{S}_a^{\alpha}, \mbf{1})^M \mbf{1}_{\alpha} \\
        &= \ab(\frac{k t_{\alpha}}{a} S^M\ab(\frac{kt_{\alpha}}{a}) \mbf{1}^{\alpha}, \mbf{1}_{\alpha}) \mbf{1}_{\alpha} \\
        &= \ab(\frac{k t_{\alpha}}{a} \1^{\alpha} , S^M\ab(\frac{kt_{\alpha}}{a})^* \mbf{1}_{\alpha}) \mbf{1}_{\alpha}  \\
        &= J^M\ab(\frac{kt_{\alpha}}{a}) \bigg|_{\pt_{\alpha}} \\
        &= I^M\ab(q, \frac{kt_{\alpha}}{a}) \bigg|_{\pt_{\alpha}}.
    \end{align*}
    Using $p|_{\pt_{\alpha}} = -t_{\alpha}$, the desired result follows from direct computation with the $I$-function.
\end{proof}

\begin{cor}\label{cor:restrictedsmat}
    The $\msc{S}_{a;i}^{\alpha}$ satisfy the following properties:
    \begin{enumerate}
        \item For any $\alpha$, we have $\msc{S}_{1;0}^{\alpha} = 1$. This implies that
            \[ \ab<\frac{\1^{\alpha}}{k t_{\alpha} - \psi}>_{0,1}^M = 0; \]
        \item For all $\alpha$, $\beta$, $a$, and $i$, the rotation property
            \[ \msc{S}_{a;i}^{\alpha} = \zeta_N^{i(\alpha-\beta)} \msc{S}_{a;i}^{\beta} \]
            holds.
        \item For all $\alpha$, $a$, and $i$, $\msc{S}_{a;i}^{\alpha}$ is a polynomial in $q$ with degree bounded by
            \[ \deg \msc{S}_{a;i}^{\alpha} \leq \begin{cases}
                \ceil{\frac{a}{k}}-1 & i < N \\
                \ceil{\frac{a}{k}} & i \geq N.
            \end{cases}
             \]
    \end{enumerate}
\end{cor}

\begin{proof}
    The first item follows by definition. The second and third items follow from the explicit form of the quantum connection given in~\Cref{lem:quantumconnection}.
\end{proof}

\subsection{Genus zero Gromov-Witten theory of $Z$}%
\label{sub:Genus zero theory of X}

We conclude this section by recalling the genus-zero Gromov-Witten theory of $Z$. Let $\tau_{Z} \coloneqq \frac{I_1}{I_0} H$ be the mirror map. Then when $\tau = \tau_{Z}$, the quantum differential equation
\[ z \d S^{Z}_{\tau}(z) = \d \tau *_{\tau} S^{Z}_{\tau}(z) \]
becomes
\[ (H+zD) S^{Z}(z)^* = S^{Z}(z)^* ( \dot{\tau}_{Z} *_{\tau_{Z}} - ) \]
using the divisor equation. Here, $\dot{\tau}_{Z} \coloneqq I_{11} H = H + D(\tau_{Z})$. Using this quantum differential equation and the mirror theorem (c.f.~\cite[Corollary 28]{applmirrortoricstack})
\[ \frac{I(q,z)}{I_0(q)} = J|_{q \mapsto qe^{\frac{I_1}{I_0}}}, \]
we can compute the quantum cohomology of $Z$ at $\tau = \tau_{Z}$, where quantum multiplication by $\tau_{Z}$ is given by the matrix
\[ A^{Z} = \begin{pmatrix}
  0 & \\
  I_{11} & 0 \\
    & I_{22} & 0 \\
  & & I_{11} & 0
\end{pmatrix}\]
in the basis $1,H,H^2,H^3$, where
\[ I_{22} = D\ab(\frac{D J_2 + J_1}{I_{11}}). \]
In particular, we find that the Yukawa coupling is given by
\[\ab<\!\ab<H,H,H>\!>_{0,3}^{Z} = p_k \frac{I_{22}}{I_{11}}. \]

The following result will be proved in~\Cref{sec:yukawa coupling}.
\begin{thm}\label{thm:yukawa}
  The normalized Yukawa coupling is given by
  \[ \frac{1}{p_k} I_0^2 I_{11}^3 \ab<\!\ab<H,H,H>\!>_{0,3}^{Z} = Y. \]
\end{thm}

Solving the ODE for the $S$-matrix, we obtain
\begin{align*}
  S^{Z}(z)^* =&{}\ I + \frac{1}{z} \begin{pmatrix}
    0 \\
    J_1' & 0 \\
           & J_2' & 0 \\
           & & J_1' & 0
  \end{pmatrix} \\ &+ \frac{1}{z^2} \begin{pmatrix}
    0 \\
    & 0 \\
    J_2 & & 0 \\
      & \frac{J_2'}{J_1'}J_1 - J_2 & &  0
  \end{pmatrix} + \frac{1}{z^3} \begin{pmatrix}
    0 \\
    & 0 \\
    & & 0 \\
    J_3 & & & 0
  \end{pmatrix},
\end{align*}
where $J_1' = I_{11}$ and $J_2' = J_1 + D J_2$.

\section{MSP theory as a CohFT}%
\label{sec:MSP theory as a CohFT}

Define the MSP $R$-matrix by the Birkhoff factorization
\begin{equation}\label{eqn:rmatrix}
    S^M(z) \begin{pmatrix}
    \on{diag}\ab\{ \Delta^{\pt_{\alpha}}(z)\}_{\alpha=1}^N & \\
    & 1
\end{pmatrix} = R(z) \begin{pmatrix}
    \on{diag}\ab\{ S^{\pt_{\alpha}}(z)\}_{\alpha=1}^N & \\
    & S^{Z}(z)
\end{pmatrix}\biggr|_{q \mapsto q'},
\end{equation}

where $\Delta^{\pt_{\alpha}}(z)$ is the Quantum Riemann-Roch~\cite{qrr} operator given by the formula
\begin{align*}
    \Delta^{\pt_{\alpha}}(z) \coloneqq \exp &\Biggl[\sum_{k \geq 0} \frac{B_{2k}}{2k(2k-1)} \Biggl(\sum_{i=1}^5 \frac{1}{(-a_i t_{\alpha})^{2k-1}}  \\
    &+ \frac{1}{(rt_{\alpha})^{2k-1}} + \sum_{\beta \neq \alpha} \frac{1}{(t_{\beta} - t_{\alpha})^{2k-1}}\Biggr) z^{2k-1} \Biggr].
\end{align*}
Here, the $B_{2k}$ are the Bernoulli numbers.

\begin{rmk}
    Let
    \[ \Delta(z) \coloneqq \begin{pmatrix}
    \on{diag}\ab\{ \Delta^{\pt_{\alpha}}(z)\}_{\alpha=1}^N & \\
    & 1
    \end{pmatrix}. \]
    We then apply~\Cref{eqn:rloc} to~\Cref{eqn:rmatrix} to obtain
    \begin{align*}
    R(z) S^{\loc}(z) &= S^M(z) \Delta(z) \\
        &= R^{\loc}(z) S^{\loc}(z) \Delta(z)  \\
        &= R^{\loc}(z) \Delta(z) S^{\loc}(z),
    \end{align*}
    which implies that $R(z) = R^{\loc}(z) \Delta(z)$.
\end{rmk}

The goal of this section is to prove the following ``theorem.''
\begin{fthm}\label{fthm:cohft}
    The MSP $[0,1]$ theory defines a CohFT 
    \[\Omega^{[0,1]}(-) \coloneqq [-]^{[0,1]} \] which is related to the local theory
    \[ \Omega^{\loc}(-) \coloneqq [-]^{\loc} \]
    by the action of the $R$-matrix defined in~\Cref{eqn:rmatrix} by
    \[ \Omega^{[0,1]} = R.\Omega^{\loc}. \]
\end{fthm}

\begin{warn}
    Because MSP theory is really a family of theories depending on a positive integer $N$ and we consider the large $N$ asymptotics of this family of theories, at the end of this section we will state a precise form of the previous ``theorem.''
\end{warn}

\subsection{Review of cohomological field theories}%
\label{sub:Review of cohomological field theories}

Cohomological field theories (CohFTs) were originally introduced by Kontsevich-Manin~\cite{km94}. The idea of an action of the symplectic loop group on CohFTs was discovered by Givental~\cite{virasorofanotoric}. Our presentation follows~\cite{relationsvia3spin}.
\begin{defn}
    A \textit{cohomological field theory} is the data $(V, \eta, \mbf{1}, (\Omega_{g,n})_{2g-2+n>0})$ of
    \begin{enumerate}
        \item A vector space $V$ with symmetric bilinear form $\eta$ and unit vector $\mbf{1} \in V$;
        \item For all $g,n$ such that $2g-2+n > 0$, a linear map
            \[ \Omega_{g,n} \colon V^{\otimes n} \to H^*(\ol{\msc{M}}_{g,n}) \]
            satisfying the following axioms:
            \begin{enumerate}
                \item $\Omega_{g,n}$ is $S_n$-equivariant, where $S_n$ acts on $V^{\otimes n}$ by permuting the tensor factors and acts on $\ol{\msc{M}}_{g,n}$ by permuting the marked points;
                \item Recall that the collection of $\ol{\msc{M}}_{g,n}$ has gluing morphisms
                    \[ q \colon \ol{\msc{M}}_{g_1, n_1+1} \times \ol{\msc{M}}_{g_2, n_2+1} \to \ol{\msc{M}}_{g_1+g_2, n_1+n_2} \]
                    which glues the two curves at their last marked points and
                    \[ s \colon \ol{\msc{M}}_{g-1, n+2} \to \ol{\msc{M}}_{g, n} \]
                    which glues the last two marked points. Let $\{e_{\alpha}\}$ be a basis for $V$. We require
                    \begin{align*}
                        s^*\Omega_{g,n}(v_1\otimes \cdots \otimes v_n) ={}& \Omega^*_{g-1, n+2}(v_1 \otimes \cdots \otimes v_n \otimes \eta^{-1}); \\
                        q^* \Omega_{g,n}(v_1 \otimes \cdots \otimes v_n) ={}& \sum_{\alpha, \beta} \eta^{\alpha\beta} \Omega_{g_1, n_1+1}(v_1 \otimes \cdots \otimes v_{n_1} \otimes e_{\alpha}) \times \\ &\times \Omega_{g_2, n_2+1}(v_{n_1+1} \otimes \cdots \otimes v_{n} \otimes e_{\beta});
                    \end{align*}
                \item Let $p \colon \ol{\msc{M}}_{g,n+1} \to \ol{\msc{M}}_{g,n}$ be the forgetful morphism. Then we require
                    \[ p^* \Omega_{g,n}(v_1 \otimes \cdots v_n) = \Omega_{g,n+1}(\mbf{1} \otimes v_1 \otimes \cdots \otimes v_n). \]
            \end{enumerate}
    \end{enumerate}
\end{defn}

There are two kinds of actions on CohFTs. The first is associated to a formal power series
\[ R(z) = 1 + R_1 z + \cdots \in \End(V) \ps{z} \]
which satisfies the symplectic condition
\[ R(z) R(-z)^* = 1. \]
We define the CohFT $R \Omega$ as a sum of stable graphs as follows. Let $G_{g,n}$ be the set of stable graphs of genus $g$ with $n$ marked points. For any $\Gamma \in G_{g,n}$, we define
\[ \on{Cont}_{\Gamma} \colon V^{\otimes n} \to H^*(\ol{M}_{g,n}) \]
as follows:
\begin{enumerate}
    \item At any vertex of $\Gamma$, we place $\Omega_{g_v, n_v}$;
    \item At the leg corresponding to the $i$-th point, we place $R^{-1}(\psi_i)$;
    \item \textit{(Flat unit).} At each edge of $\Gamma$, let $\psi'$ and $\psi''$ be the ancestor classes at the two components of the node. We insert
        \[ \frac{\eta^{-1} - R(-\psi') \eta^{-1} R^{-1}(-\psi'')^t}{\psi' + \psi''}, \]
        where the superscript $t$ denotes the matrix transpose.
\end{enumerate}
This defines a morphism $V^{\otimes n} \to H^*\ab(\prod_{v \in V(\Gamma)} \ol{\msc{M}}_{g_v, n_v})$, and we define $\on{Cont}_{\Gamma}$ to be the pushforward of this to $H^*(\ol{\msc{M}}_{g,n})$ under the natural morphism. We may finally define
\[ (R\Omega)_{g,n} \coloneqq \sum_{\Gamma \in G_{g,n}} \frac{1}{\ab|\Aut \Gamma|} \on{Cont}_{\Gamma}. \]

The second kind of action is a translation action associated to a formal power series
\[ T(z) = T_2 z^2 + T_3 z^3 + \cdots \in z^2 V\ps{z}. \]
Let $\on{st}^m \colon \ol{\msc{M}}_{g,n+m} \to \ol{\msc{M}}_{g,n}$ be the morphism forgetting the last $m$ marked points. We define
\[ (T\Omega)_{g,n}(v_1 \otimes \cdots \otimes v_n) \coloneqq \sum_{m \geq 0} \frac{1}{m!} \on{st}^m_{*} \Omega_{g,n+m} \ab(v_1 \otimes \cdots \otimes v_n \otimes \bigotimes_{j=1}^m T(\psi_{n+j})). \]

\begin{rmk}
    Note that the action of $R$ and $T$ may not preserve the flat unit (equivalently, the string equation). However, if we set
    \[ T_R = z(1-R^{-1}(z)\mbf{1})(z), \]
    then the CohFT $R.\Omega \coloneqq RT_R \Omega$ does satisfy the flat unit axiom.
\end{rmk}

\subsection{MSP $[0,1]$ theory in terms of stable graphs}%
\label{sub:MSP 01 theory in terms of stable graphs}

We will write the class
\[ [\tau_1(\psi_1), \ldots, \tau_n(\psi_n)]_{g,n}^{[0,1]} = \sum_{d \geq 0} \sum_{\Theta \in G_{g,n,d}^{[0,1]}} (-1)^{1-g} q^d \st_*^M \ab( \prod_{i=1}^n \ev_i^* \tau_i(\psi_i) \cdot \frac{[F_{\Theta}]^{\vir}}{e(N_{\Theta}^{\vir})} ) \]
as a sum of contributions from stable graphs.

Given a localization graph $\Theta$, we can stabilize it by collapsing rational bridges (called chains), rational tails (which do not carry a marked point), and ends (chains of unstable rational curves which end at a marked point). We will mark each vertex with a fixed locus $F_v$ as follows:
\begin{itemize}
    \item Every level $0$ vertex is labeled by $F_v = Z$;
    \item If $v$ has level $1$ and hour $\alpha$, then $F_v = \pt_{\alpha}$.
\end{itemize}
We will denote the resulting graph by $\Gamma = \Theta^{\st}$ and call the set of all $\Gamma$ that appear in this way by $G_{g,n}^{[0,1]}$.

\subsubsection{Tail contributions}%
\label{ssub:Tail contributions}

The tails contract into the stable vertices and contribute to the translation actions. In particular, the argument in the proof of~\Cref{lem:msponlocallagcone}, if $\psi_{\ell}$ denotes the psi class at an extra marking $\ell$ (that the tail gets contracted to), the total contribution of all possible tails attached to $\ell$ is
\[ L^{\loc}(\psi_{\ell})|_{F_v} = \psi_{\ell} + [J^M(-\psi_{\ell})]_+|_{F_v}. \]

\subsubsection{Chain and leg contributions}%
\label{ssub:Chain contributions}

Consider the two-point function
\[ W^M(z_1, z_2) \coloneqq \frac{\sum_{\alpha} e_{\alpha} \otimes e^{\alpha}}{-z_1-z_2} + \sum_{\alpha,\beta} e_{\alpha} \otimes e_{\beta} \ab<\frac{e^{\alpha}}{-z_1-\psi_1}, \frac{e^{\beta}}{-z_2-\psi_2}>_{0,2}^M. \]
By~\Cref{lem:genus0moduli} and the argument in~\cite[\S 1]{ellgwmirror} using the string equation and WDVV equation, we have
\begin{equation}\label{eqn:wtensorsmatrix}
    W^M(z_1, z_2) = - \sum_{\alpha} \frac{S^M(z_1)^{-1} e_{\alpha} \otimes S^M(z_2)^{-1} e^{\alpha}}{z_1+z_2}.
\end{equation}
We may try to compute $W^M$ by virtual localization. Every localization graph $\Theta$ that appears is a chain connecting two vertices $v_1$ and $v_2$, which carry the two insertions. At each $v_i$, there are two possibilities:
\begin{enumerate}
    \item When $v_i$ is unstable, then $\psi_i$ is a pure equivariant weight and is therefore invertible. This implies that $v_i$ will contribute terms with non-negative powers of $z_i$;
    \item When $v_i$ is stable, then $\psi_i$ is not invertible and therefore $v_i$ will contribute terms with negative powers of $z_i$.
\end{enumerate}
These contributions can be computed by taking the part of $W^M$ with either non-negative or negative powers of $z_i$. On the other hand, we may expand $W^M$ around $z_i = \infty$ by writing
\[ \frac{1}{-z_i-\psi_i} = \sum_{j\geq 0} (-1)^{j+1}\psi_i^j z_i^{-j-1}. \]
Then the coefficient of $e_{\alpha}z_i^{-j-1}$ corresponds to the correlator with insertion $e_{\alpha}\psi_i^{k}$ at $v_i$.
Using $W^M$, we may compute the chain and end contributions.

Let $e$ be an edge in $\Gamma$. We need to consider all possible chains that contract to $e$. In particular, in these contributions, the two insertions must be attached to unstable vertices, so the desired contribution is
\[ [W^M(\psi_{(e,v)}, \psi_{(e, v')})]_{+,+}. \]
Here, $v$ and $v'$ are the vertices incident to $e$. This will be inserted as a bivector at the edge $e$.

Now let $\ell$ be a leg incident to a vertex in the graph after stabilization. We need to compute the contribution of all ends with $\tau(\psi)$ as an insertion, so by the previous discussion we insert
\begin{align*}
    \Res_{z=0}&\ab([W^M(\psi_{\ell}, -z^-)]_+, \tau(z))^M \\
    &= \Res_{z=0}\ab(-\sum_{\alpha} \frac{S^M(\psi_{\ell})^{-1}e_{\alpha}\otimes S^M(-z^-)^{-1} e^{\alpha}}{\psi_{\ell-z}}\bigg\vert_+, \tau(z))^M \\
    &= [S^M(\psi_{\ell})^{-1}[( S^M(-\psi_{\ell}^-)^{-1} )^* \tau(\psi_{\ell})]_+]_+ \\
    &= [S^M(\psi_{\ell})^{-1} [S^M(\psi_{\ell}^-)\tau(\psi_{\ell})]_+]_+,
\end{align*}
where the last equality follows from the symplectic property of $S^M$.

\subsubsection{Stable graph contributions}%
\label{ssub:Stable graph contributions}

By the previous discussion, the contribution from a stable graph $\Gamma \in G_{g,n}^{[0,1]}$ is given by the following construction:
\begin{itemize}
    \item At each vertex $v \in \Gamma$, we place a linear map
        \[ \bigotimes_{j=1}^n \tau_j(z_j) \mapsto \sum_{s \geq 0} \frac{1}{s!} \st^{s}_* \ab[\prod_{j=1}^n \tau_j(\psi_j), \prod_{i=1}^{s} L^{\loc}(\psi_{n_v+i})]_{g_v, n_v+s}^{F_v, \tw}. \]
        Recall that each $L^{\loc}$ is a possible tail that contracts to $v$.
    \item At each edge of $\Gamma$, we would like to consider the bivector
        \[ [W^M(\psi_{(e,v)}, \psi_{( e,v' )})]_{+,+}, \]
        but this is not quite correct. The reason is that when we compute $W^M$ by localization, we must exclude the case where the localization graph has a single vertex. In this case, we must subtract $W^{\loc}$, which we define the same way as $W^M$ but replacing ``$M$'' with ``loc.'' Because the local theory satisfies the string and WDVV equations, the end result is
        \[ \sum_{\alpha} \frac{S^{\loc}(z)^{-1} e_{\alpha} \otimes S^{loc}(z_2)^{-1} e^{\alpha}}{z_1+z_2} - \sum_{\alpha} \frac{S^M(z)^{-1} e_{\alpha} \otimes S^M(z_2)^{-1} e^{\alpha}}{z_1+z_2} \Bigg\vert_{+,+}.\]
        Here, we view $z_1 = \psi_{(e,v)}$ and $z_2 = \psi_{(e,v')}$.
    \item At each leg, we place the contribution of all ends contracting to the leg, which was computed to be
        \[ [S^M(\psi_{\ell})^{-1} [S^M(\psi_{\ell}^-)\tau(\psi_{\ell})]_+]_+. \]
\end{itemize}

\subsection{Conversion to ancestor invariants}%
\label{sub:Conversion to ancestor invariants}

Note that while enumerative theories like Gromov-Witten theory have descendent invariants, the CohFT formalism (being a theory of cohomology classes on the moduli of stable curves) has only ancestor invariants. Recall that if $\pi^M \colon \msc{C}_{g,n,\bd} \to \msc{W}_{g,n,\bd}$ is the universal curve over the MSP moduli space and $\sigma_j \colon \msc{W}_{g,n,\bd} \to \msc{C}_{g,n,\bd}$ is the $j$-th marked point, then the descendent psi-class is given by
\[ \psi_j = c_1(\sigma_j^* \omega_{\pi^M}). \]
This class is pulled back from the moduli of \textbf{prestable} curves, which is different from the moduli of stable curves. On the other hand, if $\pi\colon \msc{C}_{g,n} \to \ol{\msc{M}}_{g,n}$ is the universal curve over the moduli of stable curves and $\sigma_j \colon \ol{\msc{M}}_{g,n} \to \msc{C}_{g,n}$ is the $j$-th marked point, the ancestor psi-class is
\[ \bar{\psi}_j = c_1(\sigma_j^* \omega_{\pi}). \]

The difference between $\psi$ and $\bar{\psi}$ comes from the process of stabilizing a prestable curves. In Gromov-Witten theory, the descendent-ancestor correspondence was found by Kontsevich-Manin~\cite{descancestor}. The version we will use is~\cite[Theorem 5.1]{virasorofanotoric}. It was proved in~\cite[Appendix 2]{qrr}, see~\cite[\S 1.5]{coatesthesis} for a more detailed exposition. Set $u = \tau^{\loc}(0,q)$ and recall that $S^{\loc}$ was an abbreviation for $S^{\loc}_u$. Because the local theory is a direct sum of twisted Gromov-Witten theories, we may directly apply the descendent-ancestor correspondence. At the vertices, we obtain
\begin{align*}
    &\sum_{s \geq 0} \frac{1}{s!} \st^s_*\ab[\prod_{j=1}^{n_v} \tau_j(z_j), L^{\loc}(\psi)^s]_{g_v, n_v+s}^{F_v, \tw} \\
    ={}& \sum_{a,b \geq 0} \frac{1}{a!b!} \st^{a+b}_*\ab[\prod_{j=1}^{n_v} \tau_j(\psi_j), ( L^{\loc}(\psi) - u )^a, u^b]_{g_v, n_v+a+b}^{F_v, \tw} \\
    ={}& \sum_{a,b \geq 0} \frac{1}{a!b!} \st^{a+b}_*\ab[\prod_{j=1}^{n_v} S^{\loc}(\bar{\psi}_j) \tau_j(\bar{\psi}_j),( S^{\loc}(\bar{\psi}) ( L^{\loc}(\bar{\psi}) - u ) )^a, u^b]_{g_v, n_v+a+b}^{F_v, \tw} .
\end{align*}
For convenience, we define
\begin{align*}
    T^{\loc}(z) &\coloneqq S^{\loc}(z) (L^{\loc}(z) - u)_+ \\
    &= S^{\loc}(z) (z(1-S^M(z)^{-1})\1|_+ - u)_+ \\
    &= z(1-R^{\loc}(z)^{-1})\1.
\end{align*}
Because $T^{\loc}(z) = O(z)$, standard facts in the intersection theory of moduli spaces of curves imply that we can replace $\bar{\psi}_1, \ldots, \bar{\psi}_{n_v}$ by the ancestor classes $\bar{\bar{\psi}}_1, \ldots, \bar{\bar{\psi}}_{n_v}$ from $\ol{\msc{M}}_{g_v, n_v}$. This implies that the vertex contribution equals
\[
    \sum_{a,b \geq 0} \frac{1}{a!b!} \st^{a+b}_*\ab[\prod_{j=1}^{n_v} S^{\loc}(\bar{\bar{\psi}}_j) \tau_j(\bar{\bar{\psi}}_j),  T^{\loc}(\bar{\psi})^a, u^b]_{g_v, n_v+a+b}^{F_v, \tw} .
\]

Using~\Cref{eqn:rloc}, we can redistribute the $S^{\loc}(z)^{-1}$ from the edges to the legs. This implies that the contribution from a stable graph $\Gamma$ is given by the following construction:
\begin{itemize}
    \item At every vertex, we insert a linear map
        \[ \bigotimes_{j=1}^{n_v} \tau_j(z_j) \mapsto \sum_{a,b \geq 0} \frac{1}{a!b!} \st^{a+b}_*\ab[\prod_{j=1}^{n_v} S^{\loc}(\bar{\bar{\psi}}_j) \tau_j(\bar{\bar{\psi}}_j),  T^{\loc}(\bar{\psi})^a, u^b]_{g_v, n_v+a+b}^{F_v, \tw}; \]
    \item At every edge we place the bivector
        \[ V^{\loc}(z,w) \coloneqq \sum_{\alpha} \frac{e_{\alpha}\otimes e^{\alpha} - R^{\loc}(z)^{-1} e_{\alpha} \otimes R^{\loc}(w)^{-1} e^{\alpha}}{z+w}; \]
    \item At every leg, we place the vector
        \[ R^{\loc}(z)^{-1} [S^M(z^-, \tau(z))]_+. \]
\end{itemize}

\subsection{Absorption of Hodge classes}%
\label{sub:Absorption of Hodge classes}

In order to formalize~\Cref{fthm:cohft}, we need to replace $R^{\loc}(z)$ with $R^{\loc}(z)$ and $[-]^{F_v, \tw}$ by $[-]^{F_v, \top}$. Following the proof of~\cite[Proposition 2.14]{relationsvia3spin}, we need to replace the translation factor $T^{\loc}(z)$ by
\begin{align*}
    z[\Delta(z)^{-1} \1 - \Delta(z)^{-1}R^{\loc}(z)^{-1}\1 + \1 - \Delta^{-1}\1] &= z (\1 - \Delta(z)^{-1}R^{\loc}(z)^{-1}\1) \\
    &= z(1-R(z)^{-1})\1 \eqqcolon T(z).
\end{align*}

We now evaluate the resulting vertex contributions, which are given by
\[ \sum_{a,b \geq 0} \frac{1}{a!b!} (\st_{a+b})_*\ab[\prod_{j=1}^{n_v} \tau_j(\bar{\bar{\psi}}_j),  T(\bar{\psi})^a, u^b]_{g_v, n_v+a+b}^{F_v, \top}. \]
When $F_v = Z$, applying the divisor equation to the insertions of $u = \frac{I_1}{I_0}q$, we obtain
\[ \sum_{a \geq 0} \frac{1}{a!b!} (\st_{a+b})_*\ab[\tau_1, \ldots, \tau_{n_v},  T(\bar{\psi})^a]_{g_v, n_v+a}^{F_v, \top}\Big\vert_{q \mapsto q e^{\frac{I_1}{I_0}}}. \]
When $F_v = \pt_{\alpha}$, the projection formula in cohomology and dimension reasons give
\begin{align*}
    \st^a_* [\tau_1, \ldots, \tau_n u^a]_{g_v, n_v+a}^{\pt_{\alpha}, \top} &= u^a f(t) \st^a_* ((\st^a)^* \tau_1 \cup \cdots \cup \tau_n) \cap [\ol{\msc{M}}_{g_v, n_v+a}] \\
    &= u^a f(t) (\tau_1 \cup \cdots \cup \tau_n) \cap \st^a_* [\ol{\msc{M}}_{g_v, n_v+a}] \\
    &= 0,
\end{align*}
where $f(t) \coloneqq \ab(\frac{1}{p_k}N(-t_{\alpha})^{N+3})^{g-1}$ is given in~\Cref{lem:topclasses}.

Let $T_v(z) \coloneqq T(z)|_{F_v}$ and set
\[ \delta_F \coloneqq \begin{cases}
    I_0(q')^{-1} & F=Z; \\
    L^{\frac{N+3}{2}}(q') & F = \pt_{\alpha}.
\end{cases}
\]
Now observe that 
\[ T_v(z) = (1-\delta_F^{-1}z) + \delta_F^{-1}z(1-\delta_F R^{-1}(z) \1) \eqqcolon (1-\delta_F^{-1}z) + \delta_F^{-1} \wt{T}_v(z), \]
where $\wt{T}_{\pt_{\alpha}}(z) = O(z^2)$ by~\Cref{eqn:r1level1} and $\wt{T}_{Z}(z) = O(z^{N-2})$ by~\Cref{eqn:r1level0}. Using the dilaton equation, we obtain
\begin{align*}
    \sum_{s \geq 0}& \frac{1}{s!} \st^s_* \ab[\prod_{i=1}^{n_v} \tau_i, T_v(\bar{\psi})^s]_{g_v, n_v+s}^{F_v, \top} \\
    &= \sum_{a,b \geq 0} \frac{1}{a!b!} \st^{a+b}_* \ab[\prod_{i=1}^{n_v} \tau_i, [1-\delta_{F_v}^{-1}\bar{\psi}]^a,[\delta_{F_v}^{-1}\wt{T}_v(\bar{\psi})]^b]_{g_v, n_v+a+b}^{F_v, \top} \\
    &= \sum_{a,b\geq 0} \frac{(1-\delta_{F_v}^{-1})^a}{b!} \binom{2g_v-2+n_v+a+b-1}{a} \st^b_* \ab[\prod_{i=1}^{n_v} \tau_i, [\delta_{F_v}^{-1} \wt{T}_v(\bar{\psi})]^b]_{g_v, n_v+b}^{F_v, \top} \\
    &= \delta_{F_v}^{2g_v-2+n_v} \sum_{b \geq 0} \st^b_* \ab[\prod_{i=1}^{n_v} \tau_i, \wt{T}_v(\bar{\psi})^b]_{g_v, n_v+b}^{F_v, \top}.
\end{align*}
Setting $N \gg 3g-3+n$, we have proved the following result:

\begin{thm}\label{thm:01cohft}
    Define the translated classes at the fixed loci by
    \begin{align*}
        [-]_{g,n}^{Z,T} &\coloneqq I_0(q')^{-(2g-2+n)} \sum_{d \geq 0} \frac{q^d e^{d\frac{I_1(q)}{I_0(q)}}}{(-t^N)^{d+1-g}} \st^Z_* (-|_Z \cap [\ol{\msc{M}}_{g,n}(Z,d)]^{\vir}); \\
        [-]_{g,n}^{\pt_{\alpha},T} &\coloneqq L(q')^{\frac{N+3}{2}(2g-2+n)} \ab(\frac{N(-t_{\alpha})^{N+3}}{p_k})^{g-1} \sum_{s \geq 0} \frac{1}{s!} \st^s_* (-|_{\pt_{\alpha}} \cup \wt{T}_{\alpha}(\psi)^k),
    \end{align*}
    where $\wt{T}_{\alpha}(z) \coloneqq z(1-L(q')^{\frac{N+3}{2}}R(z)^{-1} \1)|_{\pt_{\alpha}}$.

    When $N \gg 3g-3+n$, the $[0,1]$ classes are computed as
    \[ [\tau_1(\psi_1), \ldots, \tau_n(\psi_n)]^{[0,1]}_{g,n} = \sum_{\Gamma \in G_{g,n}^{[0,1]}} \frac{1}{\ab|\Aut\Gamma|} \on{Cont}_{\Gamma}, \]
    where $\on{Cont}_{\Gamma}$ is defined by the following construction:
    \begin{itemize}
        \item At each vertex, we place a linear map
            \[ \bigotimes_{j=1}^{n_v} \tau_j(z_j) \mapsto \ab[\prod_{j=1}^{n_v} \tau_j(\bar{\psi}_j)]_{g_v, n_v}^{F_v, T}; \]
        \item At each edge, we place the bivector
            \[ \sum_{\alpha} \frac{e_{\alpha} \otimes e^{\alpha} - R(z)^{-1} e_{\alpha} \otimes R(w)^{-1}e^{\alpha}}{z+w}; \]
        \item At each leg we place the vector
            \[ R(z)^{-1} [S^M(z^-)\tau_i(z)]_+. \]
    \end{itemize}
\end{thm}

\section{Polynomiality of the $[0,1]$ theory}%
\label{sec:Polynomiality of the 01 theory}

Our goal is to prove the following theorem, which states that the $[0,1]$ correlators are polynomials. This will be an important input to prove that the Gromov-Witten potentials are polynomial. Our strategy will be as follows:
\begin{enumerate}
    \item Prove an analogous polynomiality result for the full MSP theory;
    \item Explain how to obtain the $[0,1]$ theory from the full theory by converting localization graphs into bipartite graphs (with $[0,1]$ parts and $\infty$ parts);
    \item Prove the theorem.
\end{enumerate}

\begin{thm}\label{thm:polynomiality01}
    For any $g,n$ such that $2g-2+n > 0$ and $a_i \in [0, N+3]$, the normalized $[0,1]$ correlator
    \[ t^{\sum (n-a_i - m_i) - N(g-1)}\ab<p^{a_1} \bar{\psi}_1^{m_1}, \ldots, p^{a_n} \bar{\psi}_n^{m_n}>_{g,n}^{[0,1]} \]
    is a polynomial in $\Q[q']$ of degree at most $g-1+\frac{1}{N}\ab(3g-3+\sum a_i)$.
\end{thm}

\subsection{Polynomiality of the full theory}%
\label{sub:Polynomiality of the full theory}

\begin{lem}\label{lem:fullpolynomiality}
    Let $\tau_i \in \msc{H}$ be pure degree insertions. The full MSP generating function
    \[ \ab<p^{a_1}\bar{\psi}^{m_1}, \ldots, p^{a_n} \bar{\psi}^{m_n}>_{g,n}^M \]
    is a polynomial in $q$ of degree at most 
    \[ g-1 + \frac{3g-3 + \sum_i a_i}{N}. \]
\end{lem}

\begin{proof}
    Note that the virtual dimension of $\msc{W}_{g,n,(d,0)}$ is $N(d+1-g)+n$, so the MSP correlator vanishes when
    \[ \sum_i (a_i + m_i) < N(d+1-g)+n. \]
    On the other hand, because we are considering ancestor invariants, the invariant also vanishes when
    \[ 3g-3+n < \sum_i m_i. \]
    Subtracting the second inequality from the first, we see that the vanishing holds when
    \[ \sum_i a_i < N(d+1-g) - (3g-3), \]
    or in other words when
    \[ d > g-1 + \frac{3g-3+\sum_i a_i}{N}. \qedhere \]
\end{proof}

\subsection{Bipartite graph decomposition}%
\label{sub:Bipartite graph decomposition}

We will now decompose the full MSP theory into $[0,1]$ parts and $(1,\infty]$ parts and organize this information into bipartite graphs. We will denote the vertices of bipartite graphs as either being black (level $\infty$) or white ($[0,1]$ contributions).

\begin{defn}
    A \textit{decorated bipartite graph} is a bipartite graph with the following decorations:
    \begin{itemize}
        \item Every vertex is decorated by a genus $g_v \geq 0$ and every black vertex is decorated by $d_{\infty[v]} \in \Z$;
        \item Every edge is decorated by an integer $a_e > 0$ and an hour $\alpha_e \in [N]$;
        \item All legs are incident to white vertices (being $(1,\rho)$ insertions);
        \item Every genus zero white vertex must have $\ab|L_v| + \ab|E_v| \geq 2$.
    \end{itemize}
    We will denote the set of white vertices by $V_w$ and the set of black vertices by $V_b$. The genus and degree of a bipartite graph $\Lambda$ are given by the formulae
    \begin{align*}
        g &= h^1(\Lambda) + \sum_{v \in V} g_v \\
        d_{\infty} &= \sum_{v \in V_b} d_{\infty[v]} .
    \end{align*}
    MSP theory imposes the constraint that
    \begin{equation}\label{eqn:rhoinfty}
    d_{\infty[v]} + \frac{1}{k} \ab(2g_v-2 - \sum_{e \in E_v} (a_e-1)) \geq 0. 
    \end{equation}
    
    This geometrically corresponds to the fact that $\rho$ is nonzero on the $(1,\infty]$ part of MSP theory.
\end{defn}

We will stabilize our bipartite graphs, which induces two morphisms
\[ \cl \colon L(\Lambda) \to V(\Lambda) \]
sending each leg to the nearest vertex that survives after stabilization and
\[ c \colon L(\Lambda) \to F(\Lambda) \cup L(\Lambda) \]
sending each leg to corresponding leg (which lives on $\cl(\ell))$ after stabilization, which was originally either a flag or a leg of $\Lambda$. Define $L_v^{\circ} \coloneqq \cl^{-1}(v)$.

\begin{exm}
    We illustrate the stabilization procedure in~\Cref{fig:stab}. In this example, we have $c(\ell_1) = (v_1, e_1)$, $c(\ell_2) = (e_2, v_2)$, and $c(\ell_3) = \ell_3$.
    \begin{figure}[htpb]
    \begin{center}
    \begin{tikzpicture}[scale=1, transform shape, every edge quotes/.style = {auto,inner sep=1pt, font=\footnotesize}]
        \node[circle,fill,label=above:{$0$}] (A) at (0.5,2) {};
        \node[circle,fill,label=above:{$v_2,5$}] (B) at (1.5,2) {};
        \node[circle,draw,label=left:{$0$}] (C) at (0,0) {};
        \node[circle,draw,label=below:{$v_1, 2$}] (D) at (1,0) {};
        \node[circle,draw,label=right:{$0$}] (E) at (2,0) {};
        \node[circle,draw,label=right:{$v_3, 6$}] (F) at (3,0) {};
        \node (G) at (0,-0.7) {$\ell_1$};
        \node (H) at (2,-0.7) {$\ell_2$};
        \node (I) at (3,-0.7) {$\ell_3$};
        \draw[-] (C) -- (A) edge["$e_1$" near start] (D);
        \draw[-] (D) edge [bend left=20] (B);
        \draw[-] (B) edge [bend left=20] (D);
        \draw[-] (E) edge ["$e_2$" near start] (B);
        \draw[-] (B) -- (F);
        \draw[-] (C) -- (G);
        \draw[-] (E) -- (H);
        \draw[-] (F) -- (I);
        \node at (4,1) {$\Rightarrow$};
        \node[circle,fill,label=above:{$v_2,5$}] (J) at (5.5,2) {};
        \node[circle,draw,label=left:{$v_1,2$}] (K) at (5.5,0) {};
        \node[circle,draw,label=right:{$v_3,6$}] (L) at (7,0) {};
        \draw[-] (J) edge [bend left=20] (K);
        \draw[-] (K) edge [bend left=20] (J);
        \draw[-] (J) edge [bend left=20] (L);
        \node (M) at (5.5,-0.7) {$\ell_1$};
        \node (N) at (6.2,2) {$\ell_2$};
        \node (O) at (7, -0.7) {$\ell_3$};
        \draw[-] (K) -- (M);
        \draw[-] (J) -- (N);
        \draw[-] (L) -- (O);
    \end{tikzpicture}
    \end{center}
    \caption{Stabilization of a bipartite graph}%
    \label{fig:stab}
    \end{figure}
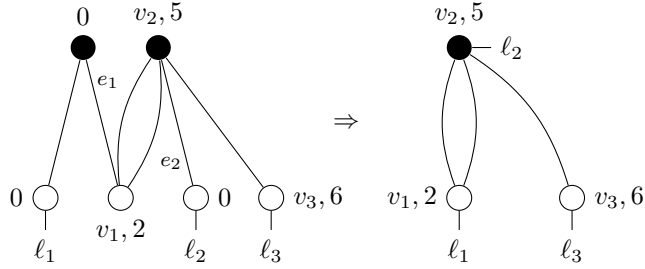
\end{exm}

We will now explain how to obtain a bipartite graph $\Lambda(\Theta)$ from a localization graph $\Theta \in G_{g,n,\bd}^{\on{reg}}$. Denote the set of all bipartite graphs obtained by this procedure by $\msc{G}_{g,n}^{d_{\infty}}$. Let $v \in V_{\infty}(\Theta)$.

\begin{defn}
    A \textit{restricted tail of $\Theta$ rooted at $v$} is the union of some $e \in E_{1\infty}$ incident to $v$ with a rational tail of $\Theta$ supported on $[0,1]$ and incident to $e$. 
\end{defn}

We first remove all restricted tails of $\Theta$ to obtain a graph $\Theta'$. We then collapse every maximal connected subgraph of $\Theta'$ supported in $[0,1]$ to a white vertex and every maximal subgraph of $\Theta'$ at level $\infty$ to a black vertex. This forms a bipartite graph $\Lambda(\Theta)$ using the induced incidence relation for edges and legs. We now assign decorations.
\begin{itemize}
    \item Each vertex is decorated by the total genus of the subgraph that was collapsed to it;
    \item The hour of an edge $e$ of $\Lambda(\Theta)$ is the hour of $e$ as an edge of $\Theta$. We also assign $a_e = -k d_e$, where $d_e = \deg \msc{L}_{\msc{C}_e}$;
    \item For every black vertex $v \in V_b$, let $\Theta_v \subset \Theta$ be the subgraph at level $\infty$ contracting to $v$. We assign
        \[ d_{\infty[v]} \coloneqq d_{\infty \Theta_v} + \sum_{\substack{e \in E_{1\infty}(\Theta) \\ (e,v) \in F(\Theta)}} d_{\infty e}. \]
        The constraint~\Cref{eqn:rhoinfty} is exactly equal to the condition that 
        \[\deg(\msc{L}^{-k} \otimes \omega_{\msc{C}}^{\log}|_{\msc{C}_{\Theta_v}}) \geq 0,\] 
        so it is always satisfied.
\end{itemize}

\begin{exm}
    We show an example of obtaining a decorated bipartite graph as follows in~\Cref{fig:bifromloc}. Note that we omit all decorations besides the genus.
    \begin{figure}[htpb]
    \begin{center}
    \begin{tikzpicture}[scale=1, transform shape, every edge quotes/.style = {auto,inner sep=1pt, font=\footnotesize}]
        \node[circle,fill,label=above:{$2$}] (A2) at (0.5,2) {};
        \node[circle,fill,label=above:{$1$}] (B2) at (1.5,2) {};
        \node[circle,fill,label=above:{$4$}] (C2) at (3,2) {};
        \node[circle,fill,label=left:{$0$}] (A1) at (0.5,0) {};
        \node[circle,fill,label=left:{$2$}] (B1) at (2,0) {};
        \node[circle,fill,label=left:{$1$}] (C1) at (4,0) {};
        \node[circle,fill,label=left:{$0$}] (A0) at (0,-2) {};
        \node[circle,fill,label=left:{$0$}] (B0) at (1,-2) {};
        \node[circle,fill,label=left:{$5$}] (C0) at (3,-2) {};
        \draw[-] (A0) -- (A1) -- (B0);
        \draw[-] (A1) -- (A2) -- (B2) -- (B1) -- (C0) edge[bend left=20] (C1);
        \draw[-] (C1) edge[bend left=20] (C0);
        \draw[-] (C0) -- (3.5,-2);
        \draw[-] (B1) -- (C2) -- (C1);
        \node at (5,1) {$\Rightarrow$};
        \node[circle,fill,label=above:{$3$}] (AB) at (6,2) {};
        \node[circle,fill,label=above:{$4$}] (BB) at (7.5,2) {};
        \node[circle,draw,label=left:{$9$}] (AW) at (7,0) {};
        \draw[-] (AB) -- (AW) edge[bend left=20] (BB);
        \draw[-] (BB) edge[bend left=20] (AW);
        \draw[-] (AW) -- (7.5,0);
    \end{tikzpicture}
    \end{center}
    \caption{Obtaining a bipartite graph from a localization graph}%
    \label{fig:bifromloc}
    \end{figure}
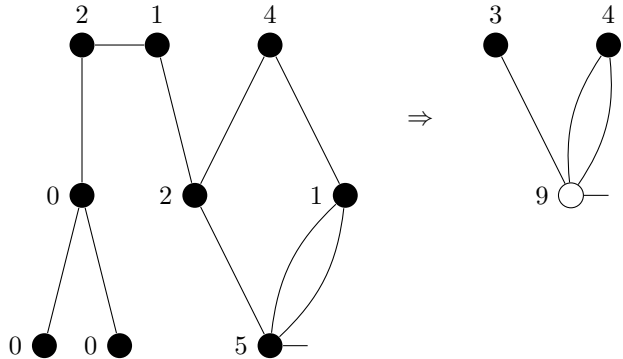
\end{exm}

\subsection{Contributions of the FJRW part}%
\label{sub:Contributions of the FJRW part}

We now analyze the contributions of the black vertices and their incident edges. Let $\Lambda \in \msc{G}_{g,n}^{d_{\infty}}$ and let $v \in V_b(\Lambda)$. Let $E_v = \{e_1, \ldots, e_m \}$. Now we let $[v]$ be the bipartite graph with one black vertex, $m$ edges $e_1, \ldots, e_m$, $m$ genus $0$ white vertices $v_1, \ldots, v_m$, and $m$ legs $\ell_1, \ldots, \ell_m$, where each $v_i$ is incident to $e_i$ and each $\ell_i$ is attached to $v_i$. We allow all $\alpha_{e_i}$ and $a_{e_i}$ to carry over. In particular, we have
\[ [v] \in \msc{G}_{g_v, m}^{d_{\infty [v]}}. \]
\begin{exm}
    In the bipartite graph in~\Cref{fig:bifromloc}, if $v$ denotes the genus $4$ black vertex, then $[v]$ is given in~\Cref{fig:vertex}. 
    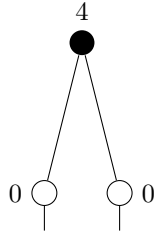
\begin{figure}[htpb]
    \begin{center}
    \begin{tikzpicture}[scale=1, transform shape, every edge quotes/.style = {auto,inner sep=1pt, font=\footnotesize}]
        \node[circle,fill,label=above:{$4$}] (AB) at (0,2) {};
        \node[circle,draw,label=left:{$0$}] (AW) at (-0.5,0) {};
        \node[circle,draw,label=right:{$0$}] (BW) at (0.5,0) {};
        \draw[-] (AW) -- (AB) -- (BW);
        \draw[-] (AW) -- (-0.5, -0.5);
        \draw[-] (BW) -- (0.5, -0.5);
    \end{tikzpicture}
    \end{center}
    \caption{$[v]$ for $v$ the genus $4$ black vertex in~\Cref{fig:bifromloc}.}%
    \label{fig:vertex}
    \end{figure}
\end{exm}

\begin{defn}
    A regular localization graph $\Theta$ \textit{strongly contracts to $[c]$} if $\Lambda(\Theta) \cong [v]$ and for each $i = 1, \ldots, m$ the subgraph $\Theta_{v_i}$ of $\Theta$ that is collapsed to $v_i$ has $d_{0\Theta_{v_i}} = 0$.
\end{defn}

\begin{defn}
    Let $G_{[v]}$ be the set of all graphs $\Theta$ strongly contracting to $[v]$. Define
    \[ \Cont_{[v]}^{\infty}\ab(\bigotimes_{i=1}^m \bar{\psi}_i^{a_i}) \coloneqq \sum_{d \geq 0} (-1)^{1-g} q^d \sum_{\substack{\Theta \in G_{[v]} \\ d_0(\Theta) = d}} \int_{[F_{\Theta}]^{\vir}} \frac{\prod_{i=1}^m \ev_i^* \1_{\alpha_i} \bar{\psi}_i^{a_i}}{e(N_{\Theta}^{\vir})}. \]
    In the case without insertions, we will abbreviate the result as $\Cont_{[v]}^{\infty}$.
\end{defn}

\begin{thm}\label{thm:decomposition}
    Let $\tau_1, \ldots, \tau_n \in \msc{H}^{\ev}$. Then the total MSP correlator is given by the formula
    \begin{align*}
        \ab<\tau_1\bar{\psi}_1^{m_1}, \ldots, \tau_n \bar{\psi}_n^{m_n}>_{g,n,d_{\infty}}^M ={}& \sum_{\Lambda \in \msc{G}_{g,n}^{d_{\infty}}} \frac{1}{\ab|\Aut \Lambda|} \prod_{v \in V_b(\Lambda)} \Cont_{[v]}^{\infty}\ab(\bigotimes_{i \in L_v^{\circ}} \bar{\psi}_{c(i)}^{m_i}) \\
        &\prod_{v \in V_w(\Lambda)} \ab<\bigotimes_{i \in L_v} \tau_i \bigotimes_{i \in L_v^{\circ}} \bar{\psi}_{c(i)}^{m_i} \bigotimes_{e \in E_v} \frac{\1^{\alpha_e}}{\frac{k t_{\alpha_e}}{a_e} - \psi_{ (e,v) }}>_{g_v, n_v}^{[0,1]},
    \end{align*}
    and $\Cont_{[v]}^{\infty}\ab(\bigotimes_{i \in L_v^{\circ}})$ is a polynomial in $q$ of degree at most
    \[ d_{\infty[v]} + \frac{1}{k}\ab(2g_v-2-\sum_{e \in E_v}(a_e-1)). \]
    Here, the unstable correlators are defined by
    \begin{align}
        \ab<\frac{x}{z_1 - \psi_1}, \frac{y}{z_i-\psi_2}>_{0,2}^{[0,1]} &\coloneqq \frac{(x,y)^M}{z_1+z_2} + \ab<\frac{x}{z_1 - \psi_1}, \frac{y}{z_2-\psi_2}>_{0,2}^{M}; \label{eqn:01twodescendents} \\
        \ab<\frac{x}{z_1 - \psi_1}, y>_{0,2}^{[0,1]} &\coloneqq (x,y)^M + \ab<\frac{x}{z_1 - \psi_1}, y>_{0,2}^{M} \label{eqn:01onedescendent} .
    \end{align}
\end{thm}

\begin{proof}
    We will give the proof first in the case in which there are no ancestor insertions. It will be clear that the argument still works with ancestor insertions. First, note that the decomposition result follows from virtual localization.

Now let $\Gamma$ be the level $\infty$ part of a regular localization graph. To be precise, there exists some $\Theta \in G^{\on{reg}}_{g,n,\bd}$ such that $\Gamma$ has vertices $V_{\infty}(\Theta)$, edges $E_{\infty\infty}(\Theta)$, and legs $E_{1\infty}(\Theta)$ with the same incidence relations as in $\Theta$. Because $\Theta$ is regular, all legs have monodromy that is FJRW-narrow.\footnote{Precisely, we mean that $m \in \{1, \ldots, k\}$ is \textit{FJRW-narrow} if $m \cdot a_i$ is not a multiple of $k$ for any $i=1, \ldots, 5$. This corresponds exactly to being a narrow insertion in FJRW theory for the Landau-Ginzburg phase of our model.} In fact, when $k=6$ we need only consider $\frac{1}{6}$ and $\frac{1}{3}$ insertions, while for $k=8,10$ we only need to consider $\frac{1}{k}$ insertions. Denoting our monodromy assignment by $\gamma = \ab(\frac{w_1}{k}, \ldots, \frac{w_m}{k})$, we see that $\Gamma \in G_{g_{\Gamma}, \gamma, \bd_{\Gamma}}^{\on{reg}}$. 

Recall that 
\[\msc{W}_{(\Gamma)} \subset \msc{W}_{g_{\Gamma}, \gamma, \bd_{\Gamma}} \] 
is the open and closed substack of $\msc{W}_{g_{\Gamma}, \gamma, \bd_{\Gamma}}^{\T}$ whose elements have localization graph $\Gamma$. We will consider the descendent classes pulled back from the coarse moduli space under the forgetful map
\[ \msc{W}_{g_{\Gamma}, \gamma, \bd_{\Gamma}} \to \mf{M}_{g,m} \]
to the moduli of prestable curves. Composing the above with $\mf{M}_{g,m} \to \ol{\msc{M}}_{g,m}$, we obtain ancestor classes. Now define the FJRW-type invariant by
\[ \msc{F}_{\Gamma}\ab(\vec{z},\bigotimes_{i=1}^m \bar{\psi}_i^{c_i}) \coloneqq q^{d_{0 \Gamma}} \int_{[\msc{W}_{(\Gamma)}]^{\vir}} \frac{1}{e(N_{\Gamma}^{\vir})} \prod_{i=1}^m \frac{\bar{\psi}_i^{c_i}}{z_i - \frac{\psi_i}{k_i}}. \]
Here, we define $k_i =k$ when $k=8, 10$ or when $\ell_i$ has monodromy $\frac{1}{6}$, and $k_i = 3$ when $\ell_i$ has monodromy $\frac{1}{3}$. Now let $\vec{a} = (a_1, \ldots, a_m)$ be a vector of non-negative integers such that $a_i \equiv w_i \pmod{5}$. Define
\[ \msc{F}_{\Gamma, \vec{a}}\ab(\bigotimes_{i=1}^m \bar{\psi}_i^{c_i}) \coloneqq (-1)^{1-g_{\Gamma}} \msc{F}_{\Gamma}\ab( \frac{-t_{\alpha_1}}{a_1}, \ldots, \frac{-t_{\alpha_m}}{a_m}; \bigotimes_{i=1}^m \bar{\psi}_i^{c_i}). \]
Again, when there are no ancestor insertions, we abbreviate this to $\msc{F}_{\Gamma, \vec{a}}$.

Now let $\Theta \in G_{[v]}$. By definition, all edges in $E_{1\infty}(\Theta)$ are one of two types:
\begin{enumerate}
    \item The first type are edges $e_1, \ldots, e_m$ such that the vertex $v_i \in V_1(\Theta)$ incident to $e_i$ is unstable;
    \item The second type are edges $e'_1, \ldots, e'_{\ell}$ corresponding to restricted tails.
\end{enumerate}

For $s=1, \ldots, m$, define $d_{\infty e_s} = \frac{a_s}{k}$ and $\alpha_{e_s} = \alpha_s$. For $i = 1, \ldots, \ell$, define $d_{\infty e'} = \frac{b_i}{k}$ and $\alpha_{e_i'} = \beta_i$. Note that when $b_j = 1$, the contributions of the $[0,1]$ parts of all possible restricted tails is
\begin{align*}
    \ab<\frac{\1^{\beta}}{k t_{\beta} - \psi}>_{0,1}^M  
    &= 0
\end{align*}
by the first item in~\Cref{cor:restrictedsmat}. Therefore, we may assume all $b_i \geq 2$.

Now note that at level $\infty$, the field $\rho$ is nowhere vanishing, so $\msc{L}^k \cong \omega_{\msc{C}}^{\log}$ for any level $\infty$ MSP field $\xi = (\msc{C}, \Sigma^{\msc{C}}, \msc{L}, \msc{N}, \varphi, \rho, \mu, \nu)$. This implies that
\[ \deg \msc{L}|_{\msc{C}_{\Gamma}} = \frac{1}{k}(2g_{\Gamma} -2 + m + \ell), \]
which translates numerically into
\begin{equation}\label{eqn:boundbs}
    0 \leq d_{0\Gamma} = d_{\infty[v]} + \frac{1}{k} (2g_{\Gamma}-2) - \frac{1}{k} \ab(\sum_{s=1}^m (a_s-1) + \sum_{i=1}^{\ell}(b_i-1)). 
\end{equation}
    Recall that $\deg \msc{F}_{\Gamma, (\vec{a}, \vec{b})} = d_{0\Gamma}$ by definition. We will bound the degree of
    \[ \Cont([v], \Gamma, \vec{b}) \coloneqq \sum_{d \geq 0} (-1)^{1-g} \sum_{\substack{\Theta \in G_{[v]} \\ \Theta_{\infty} \cong \Gamma \\ d_0(\Theta) = d \\ d_{\infty e_i'} = \frac{b_i}{5} }} \int_{[F_{\Theta}]^{\vir}} \frac{\prod_{i=1}^m \ev_i^* \1_{\alpha_i} \bar{\psi}_i^{a_i}}{e(N_{\Theta}^{\vir})}.\]
    The desired result will then follow from the fact that the sum is finite by~\Cref{eqn:boundbs} and the fact that both $d_{0\Gamma}$ and all $b_i$ are positive.

    We begin with the case where we are not in the unstable genus $0$ setting, so $(g_{\Gamma}, m+\ell, d_{0\Gamma})$ is not either $(0,2,0)$ or $(0,1,0)$. Let $A_e$ be the edge contribution to the virtual localization formula for $e \in E_{1\infty}^{\alpha}$. This is a pure equivariant weight, so its exact value is not important. Then set
    \[ c_{\vec{a}, \vec{b}}(\vec{\alpha}, \vec{\beta}) \coloneqq \frac{\ab(\prod_{s=1}^m A_{e_s}) \ab(\prod_{i=1}^{\ell} A_{e'_i})}{\ab(\prod_{s=1}^m a_s)\ab(\prod_{i=1}^{\ell} b_i)}. \]
    The virtual localization formula gives
    \begin{align*}
        \Cont([v], \Gamma, \vec{b}) ={}& \frac{c_{\vec{a}, \vec{b}}(\vec{\alpha}, \vec{\beta})}{\ab|\Aut \vec{b}|} \prod_{s=1}^m (k t_{\alpha_s}) \prod_{i=1}^{\ell}(kt_{\beta_i}) \ab<\frac{\1^{\alpha}}{\frac{kt_{\beta_i}}{b_i} - \psi}>_{0,1}^M \Bigg\vert_{\pt_{\beta_i}} \cdot 
        \msc{F}_{\Gamma, (\vec{a}, \vec{b})} \\
        ={}& \frac{c_{\vec{a}, \vec{b}}(\vec{\alpha}, \vec{\beta})}{\ab|\Aut \vec{b}|}\prod_{s=1}^m (k t_{\alpha_s}) \prod_{i=1}^{\ell} \ab((kt\beta_i)\frac{kt_{\beta_i}}{b_i} \msc{S}_{b_i;0}^{\beta_i})  
        \msc{F}_{\Gamma, (\vec{a}, \vec{b})}.
    \end{align*}
    Using the degree bound for the $\msc{S}_{b_i;0}^{\beta_i}$ from~\Cref{cor:restrictedsmat}, we obtain
    \begin{align*}
        \deg \Cont([v], \Gamma, \vec{b}) &= \deg \msc{F}_{\Gamma, (\vec{a}, \vec{b})} + \sum_{i=1}^{\ell} \deg \msc{S}_{b_i;0}^{\beta_i} \\
        &\leq d_{0\Gamma} + \sum_{i=1}^{\ell} \ab(\ceil{\frac{b_i}{k}}-1) \\
        &\leq d_{0\Gamma} + \frac{1}{k} \sum_{i=1}^{\ell} (b_i-1) \\
        &= d_{\infty[v]} + \frac{1}{k}(2g_v-2) - \frac{1}{k} \sum_{s=1}^m (a_s-1),
    \end{align*}
    where the last equality follows from~\Cref{eqn:boundbs}.

    We now consider the unstable cases. When $(m,\ell) = (0,2)$, the total degree is
    \begin{align*}
        \deg \msc{S}_{b_1;0}^{\beta_1} + \deg \msc{S}_{b_2;0}^{\beta_2} &\leq \ceil{\frac{b_1}{k}} + \ceil{\frac{b_2}{k}} - 2 \\
        &\leq \frac{b_1-1}{k} + \frac{b_2-2}{k} \\
        &= d_{\infty[v]} + \frac{1}{k}(2g_v-2).
    \end{align*}
    When $(m, \ell) = (1,1)$, the total degree is
    \begin{align*}
        \deg \msc{S}_{b;0}^{\beta} &\leq \frac{b-1}{k} \\
        &= d_{\infty[v]} + \frac{1}{k}(2g-2 - (a-1)).
    \end{align*}
    Finally, when $\ell = 0$, the contribution is a constant, and the degree bound follows from~\Cref{eqn:boundbs}.
\end{proof}

\subsection{The $[0,1]$ contributions}%
\label{sub:Proof of polynomiality of 01}

Our first step is to convert the descendent insertions at white vertices in~\Cref{thm:decomposition} into ancestors. We cannot simply use the known descendent-ancestor correspondence from Gromov-Witten theory because the MSP $[0,1]$ theory is not the Gromov-Witten theory of any target. We begin by handling the unstable contributions, which either have one leg and one flag or two flags.

For $f(z) \in \msc{H}\ps{z}$, define
\[ \mc{S} f(z) \coloneqq S^M(z^-) f(z) \in \mc{H}\ps{z, z^{-1}}. \]
This is well-defined because for each $d$ by virtual dimension constraints, only a finite number of terms of
\[ \frac{1}{z-\psi} = \sum_{i \geq 0} \psi^i z^{-i-1} \]
contribute to the coefficient of $q^d$.

\begin{lem}
    We have the formula
    \[ \ab<\frac{\tau_1}{z-\psi_1}, \tau_2(\psi_2)>_{0,2}^{[0,1]} = (S^M(z) \tau_1, \mc{S}\tau_2(-z)_+)^M. \]
\end{lem}

\begin{proof}
    Note from~\Cref{eqn:01twodescendents} that
    \begin{align*}
        \ab<\frac{x}{z_1 - \psi_1}, \frac{y}{z_i-\psi_2}>_{0,2}^{[0,1]} &= \frac{(x,y)^M}{z_1+z_2} + \ab<\frac{x}{z_1 - \psi_1}, \frac{y}{z_2-\psi_2}>_{0,2}^{M} \\
        &= \ab(x \otimes y, W^M(-z_1, -z_2)) \\
        &= \sum_{\alpha}\ab(x\otimes y, \frac{S^M(z)^* e_{\alpha} \otimes S^M(w)^* e^{\alpha}}{z+w})^M \\
        &= \sum_{\alpha} \ab(S^M(z)x \otimes S^M(w)y, \frac{e_{\alpha} \otimes e^{\alpha}}{z+w})^M \\
        &= \frac{(S^M(z) x, S^M(w)y)^M}{z+w}.
    \end{align*}
    Now we simply expand $\frac{1}{w+z} = \sum_{i\geq 0} w^{-i-1} z^i$ and compute
    \begin{align*}
        \ab<\frac{\tau_1}{z-\psi_1}, \tau_2(\psi_2)>_{0,2}^{[0,1]} &= \Res_{w=0}\ab<\frac{\tau_1}{z-\psi_1}, \frac{\tau_2(\psi_2)}{w-\psi_2}>_{0,2}^{[0,1]} \\
        &= \Res_{w=0} \frac{(S^M(z) \tau_1, \mc{S}\tau_2(w))^M}{z+w} \\
        &= (S^M(z)\tau_1, \mc{S}\tau_2(-z)_+)^M. \qedhere
    \end{align*}
\end{proof}

By the previous lemma and~\Cref{defn:specializeds}, the contribution of a genus zero white vertex with two flags is
\begin{align*}
    \ab<\frac{\1^{\alpha_1}}{\frac{k t_{\alpha_1}}{a_{e_1}} - \psi_1}, \frac{\1^{\alpha_2}}{\frac{kt_{\alpha_2}}{a_{e_2}} - \psi_2}>_{0,2}^{[0,1]} &= \frac{(\msc{S}_{a_{e_1}}^{\alpha_1}, \msc{S}_{a_{e_2}}^{\alpha_2})^M}{\frac{kt_{\alpha_1}}{a_{e_1}} + \frac{kt_{\alpha_2}}{a_{e_2}}} \eqqcolon \ab<\frac{\msc{S}_{a_{e_1}}^{\alpha_1}}{\frac{k t_{\alpha_1}}{a_{e_1}} - \bar{\psi}_1}, \frac{\msc{S}_{a_{e_2}}^{\alpha_2}}{\frac{kt_{\alpha_2}}{a_{e_2}} - \bar{\psi}_2}>_{0,2}^{[0,1]}
\end{align*}
and the contribution of a genus $0$ white vertex with one flag and one leg with insertion $\tau$ is
\begin{align*}
    \ab<\tau, \frac{\1^{\alpha}}{\frac{kt_{\alpha_2}}{a_{e_2}} - \psi_2}>_{0,2}^{[0,1]} = (\tau, \msc{S}_{a_e}^{\alpha})^M \eqqcolon\ab<\tau, \frac{\msc{S}_{a_e}^{\alpha}}{\frac{kt_{\alpha_2}}{a_{e_2}} - \bar{\psi}_2}>_{0,2}^{[0,1]}.
\end{align*}
In both cases, the degree bounds are given by~\Cref{cor:restrictedsmat}.

We now proceed to the stable case. If we consider the $[0,1]$ class
\[ [\tau_1(\psi_1), \ldots, \tau_n(\psi_n)]_{g,n}^{[0,1]}, \]
for any insertion $\tau_i(\psi_i)$, the corresponding contribution in the stable graph localization formula at the fixed locus is given by
\[ [S^M(\psi_j)^{-1} [S^M(\psi_j^-) \tau_j(\psi_j)]_+]_+ = [S^M(\psi_j)^{-1}[\mc{S}\tau_j(\psi_j)]_+]_+. \]
Applying the descendent-ancestor correspondence for the Gromov-Witten theory of the fixed locus, we obtain
\[ [S^{\loc}(\bar{\psi}_j) S^M(\bar{\psi}_j)^{-1}[\mc{S}\tau_j(\bar{\psi}_j)]_+]_+ = R^{\loc}(\bar{\psi}_j)^{-1} [\mc{S}\tau_j(\bar{\psi}_j)]_+. \]
However, an idea of Givental, for example~\cite[Theorem 9.1]{virasorofanotoric} in the semisimple case and~\cite[Theorem 3.5]{virasorotoricbundle} in the general setting, shows that localization always introduces a copy of $R^{\loc}(z)^{-1}$. Undoing this, we see that
\[ [\tau_1(\psi_1), \ldots, \tau_n(\psi_n)]_{g,n}^{[0,1]} = [\mc{S}\tau_1(\bar{\psi}_1), \ldots, \mc{S}\tau_n(\bar{\psi}_n)]_{g,n}^{[0,1]}. \]
Finally, we compute
\begin{align*}
    \ab[\mc{S}^M(z^-) \frac{x}{u-z}]_+ &= \ab[(S_0 + S_1 z^{-1} + \cdots)\frac{x}{u}\ab(1 + \frac{z}{u} + \cdots)]_+ \\
    &= \frac{1}{u} \ab(S_0x + S_1 \frac{x}{u} + \cdots) \ab(1 + \frac{z}{u} + \cdots) \\
    &= \frac{S^M(u) x}{u-z}
\end{align*}
for all $x \in \mc{H}$. Setting $x = \1^{\alpha_e}$, $u = \frac{kt_{\alpha}}{a_e}$, and $z = \bar{\psi}_{( e,v )}$, we see that the white vertex contribution from~\Cref{thm:decomposition} becomes
\[ \ab<\bigotimes_{i \in L_v} \tau_i \bigotimes_{i \in L_v^{\circ}} \bar{\psi}_{c(i)}^{m_i} \bigotimes_{e \in E_v} \frac{\msc{S}_{a_e}^{\alpha_e}}{\frac{k t_{\alpha_e}}{a_e} - \bar{\psi}_{ (e,v) }}>_{g_v, n_v}^{[0,1]}. \]

\begin{proof}[Proof of~\Cref{thm:polynomiality01}]
    We will induct on the genus $g$. Because the full MSP theory and the $[0,1]$ theory coincide when $g=0$ by~\Cref{lem:genus0moduli}, the bound holds when $g=0$.

    Now we assume that the bound holds for
    \[ \ab<p^{a_1} \bar{\psi}_1^{m_1}, \ldots, p^{a_n} \bar{\psi}_n^{m_n}>_{h,s}^{[0,1]} \]
    for any $h < g$ and any number of marked points $s$. We now consider the $[0,1]$ correlator of genus $g$ with $n$ marked points. Note that $\msc{G}_{g,n}$ contains a distinguished graph $\Lambda_{g,n}$ with a single white vertex of genus $g$ and $n$ marked point. Set $\msc{G}_{g,n}^{\circ} \coloneqq \msc{G}_{g,n} \setminus \{\Lambda_{g,n}\}$. Then~\Cref{thm:decomposition} implies that
    \[ \ab<p^{a_1} \bar{\psi}_1^{m_1}, \ldots, p^{a_n} \bar{\psi}_n^{m_n}>_{h,s}^{[0,1]} = \ab<p^{a_1} \bar{\psi}_1^{m_1}, \ldots, p^{a_n} \bar{\psi}_n^{m_n}>_{h,s}^{M} - \sum_{\Lambda \in \msc{G}_{g,n}^{\circ}} \Cont_{\Lambda}. \]
    The first term on the right hand side satisfies the degree bound by~\Cref{lem:fullpolynomiality}. It remains to prove the degree bound for the second term on the right hand side. Fortunately, this is a relatively simple degree-counting argument.

    Let $\Lambda \in \msc{G}_{g,n}^{\circ}$. By our previous discusssion, the contribution from $\Lambda$ is (up to a constant factor)
    \begin{equation}\label{eqn:bipartitecontancestor}
    \prod_{v \in V_b(\Lambda)} \Cont_{[v]}^{\infty}\ab(\bigotimes_{i \in L_v^{\circ}} \bar{\psi}_{c(i)}^{m_i}) 
        \prod_{v \in V_w(\Lambda)} \ab<\bigotimes_{i \in L_v} p^{a_i} \bigotimes_{i \in L_v^{\circ}} \bar{\psi}_{c(i)}^{m_i} \bigotimes_{e \in E_v} \frac{\msc{S}_{a_e}^{\alpha_e}}{\frac{k t_{\alpha_e}}{a_e} - \bar{\psi}_{ (e,v) }}>_{g_v, n_v}^{[0,1]} .
    \end{equation}
        We will assume this is nonzero. Recall that 
        \[ \mc{S}_{a_e}^{\alpha_e} = \sum_{j=0}^{N+3} \msc{S}_{a_e; j}^{\alpha_e} \phi^j. \]
        We will now count the degree of~\Cref{eqn:bipartitecontancestor}.
        \begin{itemize}
            \item Suppose $v \in V_w(\Lambda)$ is a stable white vertex. By the proof of~\Cref{lem:genus0moduli}, the total genus of all black vertices in $\Lambda$ is positive, so $g_v < g$. Therefore, we may apply the induction hypothesis. Because
                \[ \phi^j = \frac{1}{p_k}(p^{N+3-j} + \cdots), \]
                we see that
                \begin{align*}
                    \deg \msc{S}_{a_e; j}^{\alpha_e} + \frac{1}{N} \deg \phi^j &\leq \ceil{\frac{a_e}{k}} + \frac{N+3-j}{N} \\
                    &\leq \frac{a_e-1}{k} + 1 + \frac{3}{N}
                \end{align*}
                by~\Cref{cor:restrictedsmat}. Therefore, the total contribution at $v$ has degree at most
                \[ g_v - 1 + \sum_{e \in E_v} \frac{a_e-1}{k} + \ab|E_v| + \frac{3g_v-3 + 3\ab|E_v|+\sum_{i \in L_v} a_i}{N}. \]
            \item Suppose $v \in V_w(\Lambda)$ is unstable. Then
                \begin{itemize}
                    \item If $v$ has one edge $e$ and one insertion $p^{a_i} \bar{\psi}_i^{m_i}$, the $[0,1]$ correlator is a polynomial of degree at most
                        \begin{align*}
                            \ceil{\frac{a_e}{k}} - \delta_{a_i < N} &= \ceil{\frac{a_e}{k}} - 1 + \floor{\frac{a_i}{N}} \\
                            &\leq \frac{a_e-1}{k} + \frac{a_i}{N};
                        \end{align*}
                    \item If $v$ has two edges, then the precise statement of~\Cref{cor:restrictedsmat} implies (by degree reasons) that the unstable $[0,1]$ correlator has degree at most
                        \[ \ceil{\frac{a_{e_1}}{k}} + \ceil{\frac{a_{e_2}}{k}} - 1 \leq \frac{a_{e_1}-1}{k} + \frac{a_{e_2}-1}{k} + 1 + \frac{3}{N}. \]
                \end{itemize}
            \item Suppose that $v \in V_b(\Lambda)$ is a black vertex. By~\Cref{lem:genus0moduli}, we know $g_v \geq 1$. By~\Cref{thm:decomposition}, we know that
                \[ \deg \Cont_{[v]}^{\infty} \leq d_{\infty[v]} + \frac{1}{k}\ab(2g_v-2-\sum_{e \in E_v}(a_e-1)). \]
        \end{itemize}
        Recalling the abbreviation $\ab<->_{g,n}^{M} = \ab<->_{g,n,0}^M$ (or in other words that the total $d_{\infty}$ is zero), we see that
        \[ \sum_{v \in V_b(\Lambda)} d_{\infty[v]} = 0. \]
        Therefore, the total degree of $\Cont_{\Lambda}$ is at most
        \begin{align*}
            &\sum_{v \in V_w(\Lambda)} (g_v-1 + \ab|E_v|) + \frac{\sum_{v \in V_w(\Lambda)} 3(g_v-1 + \ab|E_v|) + \sum_{i=1}^n a_i}{N} \\
            &\qquad + \sum_{v \in V_b(\Lambda)} \frac{2}{k}(g_v-1) \\
            \leq {}& \sum_{v \in V_w(\Lambda)} (g_v-1) + \ab|E(\Lambda)| + \sum_{v\in V_b(\Lambda)} (g_v-1) \\
            &\qquad + \frac{\sum_{v \in V(\Lambda)} 3(g_v-1) + 3\ab|E(\Lambda)| + \sum_{i=1}^n a_i}{N} \\
            ={}& g-1 + \frac{3(g-1) + \sum_{i=1}^n a_i}{N}. \qedhere
        \end{align*}
\end{proof}

\section{$R$-matrix computations}%
\label{sec:Computation of the R-matrix}

This section collects various computations with the $R$-matrix that are needed in the other parts of this paper. For simplicity, we will specialize $t$ such that $t^N = -1$.

\subsection{Level $0$}%
\label{sub:The R matrix at level 0}

This subsection computes the formulae at level $0$ needed for~\Cref{sec:MSP theory as a CohFT}.

We begin by computing the tail contribution at level $0$. This is given by
\begin{align*}
    L^{\loc}(z)|_Z ={}& z \1 + J^M(0,-z)|_{Z,+} \\
    ={}& - \sum_{d=1}^{\infty} z q^d \frac{\prod_{m=1}^{kd}(kH-mz)}{\prod_{i=1}^5 \prod_{m=1}^{a_i d}(a_i H - mz) \prod_{m=1}^d ((H-mz)^N - t^N)} \Bigg\vert_+ \\
    ={}& -\sum_{d=1}^{\infty} z (q')^d\frac{\prod_{m=1}^{kd}(kH-mz)}{\prod_{i=1}^5 \prod_{m=1}^{a_i d}(a_i H - mz)} \\
    & \cdot \ab(1 + O(z^N) + H O(z^{N-1}) + H^2 O(z^{N-2}) + H^3 O(z^{N-3})) \big\vert_+ \\
    ={}& z + I^Z(q',-z)|_+ + O(z^{N-2}) \\
    ={}& z(1-I_0(q')) + I_1(q')H + O(z^{N-2}),
\end{align*}
where we have used~\Cref{thm:mirror} for the explicit formula for $J^M$, expanded $((H-mz)^N-t^N)^{-1}$ using the fact that $t$ is a unit, and finally used the fact that $H^4 = 0$.

We also compute
\begin{align*}
    \tau^{\loc}(0,q)|_Z &= \sum_{\alpha, n} \frac{1}{n!} e^{\alpha} \ab<e_{\alpha}, \1, L^{\loc}(-\psi)|_Z^n>_{0,n+2}^{Z,\tw} \\
    &= \sum_{\alpha, n} \frac{1}{n!} e^{\alpha} \ab<e_{\alpha}, \1, ((1-I_0(q'))\psi + I_1(q')H)^n>_{0,n+2}^{Z,\tw} \\
    &= \sum_{a,b,\alpha} \frac{1}{a!b!} e^{\alpha} \ab<e_{\alpha}, \1, ( (1-I_0(q'))\psi )^a, (I_1(q')H)^b>_{0,a+b+2}^{Z,\tw} \\
    &= \sum_{a,b,\alpha} \frac{(1-I_0(q'))^a}{a!b!} \binom{a+b-1}{a} e^{\alpha} \ab<e_{\alpha}, \1, (I_1(q')H)^b>_{0,a+b+2}^{Z,\tw} \\
    &= \sum_{\alpha, b} \frac{1}{b!} e^{\alpha} \ab<e_{\alpha}, 1, \tau_Z(q')^b>_{0,k+2}^{Z,\tw} \\
    &= \tau_Z(q'),
\end{align*}
where the second-to-last equality follows from the formal power series for $\frac{1}{(1-x)^n}$ and the last equality follows from the string equation.

We will now use the quantum connection~\Cref{lem:quantumconnection} and~\Cref{eqn:rmatrix} defining the $R$-matrix to compute $S^M$ and $R$ at level $0$. We first note that
\begin{align*}
    S^M(z)^*\1|_{Z,+} &= z^{-1} J^M(z)|_{Z,+} \\
    &= I_0(q') + O(z^{N-3}).
\end{align*}
Applying the quantum connection, we obtain
\begin{align*}
    S^M(z)^*p|_{Z,+} &= (D^p S^M(z)^* \1)|_{Z,+} \\
    &= HI_0(q') + zD(I_0(q')) + D(I_1(q'))H + O(z^{N-2})
\end{align*}
from the explicit computations at the beginning of this subsection. Applying~\Cref{eqn:rmatrix} and the symplectic property, we formally obtain
\[ R(z)^{-1}x|_Z = S^Z(q',z)(S^M(z)^{-1}x)|_Z \]
for all $x \in \msc{H}$. This implies that
\begin{align}
    R(z)^*\1|_Z &= I_0(q') + O(z^{N-3}); \label{eqn:r1level0} \\
    R(z)^*p|_Z &= z D(I_0(q')) + I_0 I_{11}H + O(z^{N-2}). \label{eqn:rplevel0} 
\end{align}

\begin{rmk}
    The estimate of $O(z^{N-2})$ in the formula for $S^M(z)^*p|_{Z,+}$ and the estimates of $O(z^{N-3})$ and $O(z^{N-2})$ in~\Cref{eqn:r1level0} and~\Cref{eqn:rplevel0}, respectively, come from more careful study of the formula for $J^M$ appearing in the intermediate computations of $L^{\loc}(z)$.
\end{rmk}

\subsection{Level $1$}%
\label{sub:Level 1}

In this subsection, we compute enough of the $R$-matrix at level $1$ to obtain the explicit formulae in~\Cref{sec:MSP theory as a CohFT}.

We will prove the following key facts.
\begin{prop}\label{prop:rmatlevel1}
    Let $L \coloneqq (1+rx)^{\frac{1}{N}}$. Then
    \begin{enumerate}
        \item We have the formula
    \[ \tau^{\loc}(0,q) = \tau_{\alpha}(q') = -t_{\alpha}\int_0^{q'} (L(x)-1) \frac{\d x}{x}; \]
\item The quantity \[ (R_m)_0^{\alpha} \coloneqq L^{\frac{N+3}{2}} L_{\alpha}^m (R_m \mbf{1}^{\alpha}, \mbf{1})^M \]
    is independent of $\alpha$ and polynomial in $Y$ of degree at most $k$;
        \item The tail contribution at $\pt_{\alpha}$ is given by
            \begin{align*}
                T(z)|_{\pt_{\alpha}} ={}& \mbf{1}_{\alpha} z - L^{- \frac{N+3}{2}} \mbf{1}_{\alpha}z \\
                & \cdot \ab[1 - \ab(\ab(\frac{N}{24} + C_k) + \frac{Y-1}{N} \ab(\frac{N^2}{12} + \frac{23N}{24} + \frac{47}{24})) ], 
            \end{align*}
            where $C_6 = \frac{23}{72}$, $C_8 = \frac{29}{96}$, and $C_{10} = \frac{31}{120}$.
    \end{enumerate}
\end{prop}

\begin{defn}
    Define $m \in \ab\{1, \ldots, k\}$ to be \textit{ordinary} if either $m$ is FJRW-narrow or if $m=k$.
\end{defn}

First, the MSP $I$-function satisfies the equation\footnote{From the perspective of computations, we use this equation instead of the more natural-looking equation
\[  \prod_{i=1}^5 \prod_{j=0}^{a_i-1} (a_i D_p-jz) \prod_{\alpha=1}^N(D_p+t_{\alpha}) -  q \prod_{j =1}^k (kD_p + jz)) \]
in order for the equation to have the same computational properties as the Picard-Fuchs equation for the MSP theory of the quintic in~\cite[\S 6]{nmsp2}.}
\[ D_p^5 \prod_{\alpha=1}^N(D_p+t_{\alpha}) - \frac{k^{k-5}}{\prod_i a_i^{a_i}} q \prod_{j \text{ ordinary}} (kD_p + jz)). \]
By the mirror theorem and~\Cref{eqn:rmatrix}, we obtain
\begin{equation}\label{eqn:pfrmatrix}
    D_{L_{\alpha}}^5 \prod_{\alpha=1}^N(D_{L_{\alpha}}+t_{\alpha}) - q \prod_{j\text{ ordinary}} (kD_{L_{\alpha}} + mz)) R^*(z) \1 |_{\pt_{\alpha}},
\end{equation}
where $L_{\alpha} \coloneqq -t_{\alpha} + D \tau_{\alpha}$ and $D_{L_{\alpha}} = zD + L_{\alpha}$. The $z^0$-coefficient of~\Cref{eqn:pfrmatrix} is
\[ \ab( \prod_{i=1}^5 a_i^{a_i} )L_{\alpha}^k (L_{\alpha}^N - t^N) - q \cdot k^k L_{\alpha}^k = 0. \]
The nonzero solutions are
\[ L_{\alpha} = \zeta_N^{\alpha}t \ab(1+\frac{rq}{t^N})^{\frac{1}{N}} \]
with multiplicity $1$ each.
Using the initial condition $\tau_{\alpha}|_{q=0} = 0$, we obtain~\Cref{prop:rmatlevel1} (1).

\subsubsection{Expanding the Picard-Fuchs equation}%
\label{ssub:Expanding the Picard-Fuchs equation}

Continuing our study of~\Cref{eqn:pfrmatrix}, the $z^1$-coefficient is given by
\[ \frac{N+3}{2} rq R_0^* \1 |_{\pt_{\alpha}} + N (rq + t^N)D(R_0^* \1 |_{\pt_{\alpha}}) = 0. \]
Using separation of variables, we obtain
\begin{equation}\label{eqn:r1level1}
R_0^* \1|_{\pt_{\alpha}} = L^{-\frac{N+3}{2}} 
\end{equation}

using the initial condition that
\[ R_0^*\1|_{\pt_{\alpha},q=0} = [z^0] \Delta^{\pt_{\alpha}}(z) = 1. \]
We may repeat this strategy by considering the coefficients of $z^m$ for any $m \geq 2$ and solve $R_m^*\1|_{\pt_{\alpha}}$ up to a constant. We will use the fact that $R$ is the asymptotic expansion of an oscillating integral to fix the initial conditions. 

Rewrite
\[ R(z)^* \1|_{\pt_{\alpha}} = L^{-\frac{N+3}{2}} \ab(1 + \frac{r_1}{L_{\alpha}}z + \frac{r_2}{L_{\alpha}^2} z^2 + \cdots). \]
If we denote the differential operator annihilating $R(z)^* \1 |_{\pt_{\alpha}}$ by $L_{\alpha}^{N+5} \cdot \msc{PF}$, then~\Cref{eqn:pfrmatrix} becomes
\begin{equation}
    L_{\alpha}^{N+5} \cdot \msc{PF} \ab(L^{-\frac{N+3}{2}} \cdot \ab(1 + \frac{r_1}{L_{\alpha}} z + \frac{r_2}{L_{\alpha}} z^2 + \cdots)) = 0.
\end{equation}
Setting $X = 1- L^{-N}$, we see that
\begin{align*}
  \msc{PF} &= L_{\alpha}^{-5} \ab(\frac{1}{L_{\alpha}^N} D_{L_{\alpha}}^N - (1-X))D_{L_{\alpha}}^5 + \frac{X}{r L_{\alpha}^5} \prod_{j\text{ ordinary}} (k D_{L_{\alpha}} + mz) \\
  &= \sum_{m=1}^{N+5} \msc{PF}_m z^m,
\end{align*}
where
\begin{equation}
  \msc{PF}_m = \sum_{j=0}^m f_{m,j}(X)D^{m-j}, \qquad f_{m,j} \in \Q[x]_j.
\end{equation}
For example, we have
\begin{align*}
    \msc{PF}_1 ={}& ND + \frac{N+3}{2} X  \\
    \msc{PF}_2 ={}& \frac{N(N+9)}{2}D + \frac{N^2+12N+23}{2}XD \\
    &+ \ab(\frac{N^2}{12} + \frac{25N}{24} + \frac{35}{12} + \frac{47}{24N})X^2 - \ab(\frac{N^2}{12} + N + c_k) X,
\end{align*}
where $c_6 = \frac{41}{18}$, $c_8 = \frac{217}{96}$, and $c_{10} = \frac{133}{60}$.

\begin{rmk}\label{rmk:derivatives}
    In the above, we use several important observations:
    \begin{enumerate}
        \item We have $DL_{\alpha} = \frac{L_{\alpha}X}{N}$ and $DX = X(1-X)$, so commuting $X$ or $L_{\alpha}$ to the left of the derivative increases the $X$-degree by $1$;
        \item In each monomial of $\msc{PF}_m$, the term $D_{L_{\alpha}}$ appears at most $m$ times;
        \item We have
            \begin{align*}
                D_{L_{\alpha}}^n &= L_{\alpha}^n + z \ab(nL_{\alpha}^{n-1}D + \frac{n(n-1)}{2} \frac{X}{N}) + \cdots \\
                &= \sum_{m=0}^n\ab( z^j L_{\alpha}^{n-m} \sum_{j=0}^m c_{m,j}(X)D^{m-j}),
            \end{align*}
            where $c_{m,0} = \binom{n}{m}$ and $c_{m,j} \in XQ[X]_{j-1}$ whenever $j > 0$.
    \end{enumerate}
\end{rmk}

Now, we set
\begin{align*}
  \wt{\msc{PF}} &= L^{\frac{N+3}{2}} \cdot \msc{PF}\cdot L^{-\frac{N+3}{2}} \\
  &= \msc{PF}|_{D \mapsto D - \frac{N+3}{2N} X}.
\end{align*}
Expanding $\wt{\msc{PF}}$ in terms of $z$, we see for example that $\wt{\msc{PF}}_1 = ND$. The Picard-Fuchs equation now becomes
\[ \wt{\msc{PF}} \ab(1 + \frac{r_1}{L_{\alpha}}z + \frac{r_2}{L_{\alpha}^2} z^2 + \cdots) = 0. \]
We now obtain
\[ \frac{1}{L_{\alpha}} \wt{\PF}_1 \ab(\frac{r_m}{L_{\alpha}^m}) + \frac{1}{L_{\alpha}^2} \wt{\PF}_2\ab(\frac{r_{m-1}}{L_{\alpha}^{m-1}}) + \cdots + \frac{1}{L_{\alpha}^{m+1}} \wt{\PF}_{m+1}(1) = 0, \]
which gives us
\begin{align*}
  ND \ab(\frac{r_m}{L_{\alpha}^m}) &= - \sum_{j=0}^{m-1} \frac{1}{L_{\alpha}^{m+1-j}} \wt{\PF}_{m+1-j} \ab(\frac{r_j}{L_{\alpha}^j}) \\
                              &= -\frac{1}{L_{\alpha}^m} \sum_{j=0}^{m-1} \msc{PF}_{m+1-j}|_{D \mapsto D-\frac{j}{N}X}(r_j) \\
  &\eqqcolon \frac{1}{L_{\alpha}^m}\msc{P}_m(r_0, \ldots, r_{m-1}),
\end{align*}
where we set $r_0 = 1$.

This implies that we may compute $r_m$ via the integral
\[ \frac{r_m}{L_{\alpha}^m} = \int \frac{1}{L_{\alpha}^k} \msc{P}_m(r_0, \ldots, r_{m-1}) \frac{\d L_{\alpha}}{X L_{\alpha}}. \]
This is well-defined because by~\cite[Equation 6.11]{nmsp2}, $\msc{P}_m(r_0, \ldots, r_{m-1})$ is divisible by $X$. To prove independence in $\alpha$, we will prove that the integration constants vanish.

\begin{lem}\label{lem:vanishing}
    The coefficient of $z^m$ in the series $L^{\frac{N+3}{2}} R(z)^* \mbf{1}_{\pt_{\alpha}}$ approaches $0$ when $L_{\alpha}^{-1} \to 0$.
\end{lem}

\subsubsection{Asymptotic expansion of oscillating integral}%
\label{ssub:Asymptotic expansion of oscillating integral}

Recall that the genus-$0$ MSP theory is equivalent to the equivariant GW theory of a degree $k$ hypersurface in $\P \coloneqq \P(\mbf{a}, 1^N)$, or equivalently, the $\msc{O}(k)$-twisted GW theory of $\P$. Here, we will always assume that $k < \sum_i a_i$, or equivalently that we are in the Fano case. Denote the $S$-matrix of this theory by $S^{\lambda}(z)$ and then define
\[ \Delta^{\lambda,\alpha}(z) = \exp \sum_{k>0} \frac{B_{2k}}{2k(2k-1)} \ab(\frac{1}{( k\lambda_{\alpha} )^{2k-1}} + \sum_{\beta \neq \alpha} \frac{1}{(a_{\beta} (\lambda_{\beta} - \lambda_{\alpha}))^{2k-1}})z^{2k-1}  \]
and $\Delta^{\lambda}(z) = \on{diag} \ab\{ \Delta^{\lambda,\alpha}(z)\}_{\alpha=0}^N$. Note that we will only need to consider the untwisted sector of this theory in this paper.

Following~\cite{htmstoric}, consider the Landau-Ginzburg potential $W \colon (\C^{\times})^{n+1} \to \C$ for the total space of $\msc{O}(-k)$ on $\P(a_0, \ldots, a_n)$:
\[ W(x_0, \ldots, x_n) \coloneqq \sum_{i=0}^n (x_i - \lambda_i \log x_i) - y, \qquad y^k q = \prod_{i=0}^n x_i^{a_i}. \]
To find the critical points $\mbf{x}_{\alpha}$ of $W$, we note that
\[ x_i \pdv{W}{x_i} = x_i - \lambda_i - a_i \frac{y}{k} = 0 \]
and set $\frac{y}{k} \eqqcolon L_{\alpha}$. We then obtain the equation
\[ \prod_i (a_i L_{\alpha}+\lambda_i)^{a_i} = (kL_{\alpha})^k q \]
and the critical points
\[ (x_i)_{\alpha} = a_i L_{\alpha} + \lambda_i, \qquad y_{\alpha} = k L_{\alpha}. \]
Denote the corresponding critical value by $u_{\alpha}$.

Near each critical point, we have a Lefschetz thimble $\gamma_{\alpha}$, which is a real $n$-dimensional cycle in $(\C^{\times})^{n+1}$ that becomes the vanishing cycle when restricted to the Milnor fiber. By~\cite[Lemma 6.2]{htmstoric}, there is a bijection between critical points and pairs $(\sigma, \chi)$, where $\sigma$ is a maximal cone of the fan defining $\P(a_0, \ldots, a_n)$ and $\chi \colon N(\sigma) \to \C^{\times}$ is a character of the stabilizer group $N(\sigma)$ of the torus-fixed point corresponding to $\sigma$. We then consider the $z \to 0^-$ (along the negative real axis) asymptotics of the oscillating integral
\[ \msc{I}_{\alpha}(q,z) \coloneqq \int_{\gamma_{\alpha}} e^{\frac{W}{z}} \frac{\d x_0 \wedge \cdots \d x_n}{x_0 \cdots x_n}. \]
For critical points corresponding to $\sigma$ with $N(\sigma) = 1$, or equivalently $i$ such that $a_i = 1$,~\cite[Proposition 6.9]{htmstoric} implies that
\[ (-kzD) \msc{I}_{\alpha}(q,z) \asymp e^{\frac{u_{\alpha}}{z}} (-2\pi z)^{\frac{n+1}{2}} R^{\lambda}(z)^* \1|_{\pt_{\alpha}}, \]
where $R^{\lambda}(z)$ is defined by the Birkhoff factorization
\begin{equation}\label{eqn:rmattwisted} S^{\lambda}(z) \Delta^{\lambda}(z) = R^{\lambda}(z) e^{\frac{u}{z}}. \end{equation}
Here, $u = \on{diag} \ab\{u_{\alpha}\}_{\alpha=0}^n$ in the fixed point basis $\1_{\pt_{\alpha}}$ (here, we will not need to consider the twisted sectors of $\P$).

We will now compute the asymptotic expansion using the saddle point method. Near $\mbf{x}_{\alpha}$, $W$ has the form
\[ W = u_{\alpha} + \frac{1}{2} \msc{Q}(\xi) + \sum_{k \geq 3} \frac{1}{k!} \sum_{i+1, \ldots, i_k} \partial_{i_1} \cdots \partial_{i_k} W(\mbf{x}_{\alpha}) \xi_{i_1} \cdots \xi_{i_j}, \]
where $\xi = \mbf{x} - \mbf{x}_{\alpha}$ are the local coordinates and $\msc{Q}(\xi)$ is the quadratic form given by the Hessian of $W$ at $\mbf{x}_{\alpha}$. Writing $\partial_I \coloneqq \partial_{i_1} \cdots \partial_{i_n}$ and $\xi_I \coloneqq \xi_{i_1} \cdots \xi_{i_n}$ for $I = (i_1, \ldots, i_k)$, we see that
\begin{align*}
    \msc{I}_{\alpha}(q,z) &= \int_{\gamma_{\alpha}} e^{\frac{u_{\alpha}}{z}} e^{\frac{1}{2} \frac{1}{z} \msc{Q}(\xi)} e^{\frac{1}{z}\sum_{k \geq 3} \frac{1}{k!} \sum_{\ab|I| = k} \partial_I W(\mbf{x}_{\alpha}) \xi_I} \frac{\d \xi_0 \cdots \d\xi_0}{\prod_{i=0}^n ((x_i)_{\alpha} + \xi_i)} \\
    &\asymp e^{\frac{u_{\alpha}}{z}}(-2\pi z)^{\frac{n+1}{2}} \Psi \cdot \ab(1 + \sum_{\ell > 0} f_{\ell}(-z)^{\ell}),
\end{align*}
where we set
\[ (2\pi)^{\frac{n+1}{2}}\Psi \coloneqq \prod_i \frac{1}{( x_i )_{\alpha}} \int_{\R^{n+1}} e^{-\frac{1}{2} \msc{Q}(\xi)} \d \xi_1 \cdots \d\xi_n. \]
Note that $\frac{1}{(x_i)_{\alpha} + \xi_i}$ is expanded as $\sum_{m\geq 0} (-1)^m \frac{\xi_i^m}{(x_i)_{\alpha}^{m+1}}$.

\begin{lem}
    Let $m = m_0 + \cdots + m_n$. If $m$ is even, then
    \[ L_{\alpha}^{-\frac{n+1+m}{2}} \int_{\R^n} \xi_0^{m_0} \cdots \xi_n^{k_n} e^{\frac{1}{2z} \msc{Q}(z)} \d \xi_0 \cdots \d \xi_n = (-z)^{\frac{n+1+m}{2}} \]
    as $L_{\alpha}^{-1} \to 0$. When $m$ is odd, the integral vanishes.
\end{lem}

\begin{proof}
    At $\mbf{x}_{\alpha}$ we have
    \[ \pdv{W}{x_i,x_j} = \frac{\delta_{ij}}{(a_i L_{\alpha}+\lambda_i)^2} - \frac{a_i a_j L_{\alpha}}{k (a_i L_{\alpha} + \lambda_i)(a_j L_{\alpha} + \lambda_j)}. \]
    This yields
    \begin{align*}
        \det \msc{Q} &= \prod_i \frac{1}{a_iL_{\alpha}+\lambda_i} \ab(1- \frac{L_{\alpha}}{k} \sum_i \frac{a_i}{L_{\alpha} + \frac{\lambda_i}{a_i}}) \\
        &= \ab( \prod_i \frac{1}{a_i} ) L_{\alpha}^{-(n+1)} \ab(1-\frac{\sum_i a_i}{k} + O(L_{\alpha}^{-1})).
    \end{align*}
    The rest of the proof proceeds exactly as in~\cite[Lemma 6.4]{nmsp2}.
\end{proof}

\begin{lem}\label{lem:asymptotics}
    For all $\ell > 0$, $f_{\ell}$ is a rational function of $L_{\alpha}$ satisfying
    \[ \lim_{L_{\alpha}^{-1} \to 0} f_{\ell}(L_{\alpha}) = 0. \]
\end{lem}

\begin{proof}
    The proof is the same as~\cite[Proposition 6.1]{nmsp2}.
\end{proof}

\subsubsection{Specialization}%
\label{ssub:Specialization}

We now set $n+1 = N+5$ and specialize the equivariant parameters to $0,0,0,0,0,-\zeta_N t, -\zeta_N^2 t, \ldots, -\zeta_N^N t$. Then we obtain
\[ L_{\alpha}^k(L_{\alpha}^N - t^N) = r q L_{\alpha}^k. \]
Therefore, we have $N$ critical points corresponding to
\[ L_{\alpha} = \zeta_N^{\alpha}(t^N + rq)^{\frac{1}{N}}. \]
These correspond to the fixed points at level $1$ in the MSP theory, so we may use the arguments in the previous subsection. In particular, we have
\[ \det \msc{Q} = \frac{-N}{k} \prod_i \frac{a_i}{(a_i L_{\alpha} + \lambda_i)^2} L_{\alpha}^{\frac{N+5}{2}} \]
and therefore
\begin{align*}
    \Psi &= \prod_i \frac{1}{(x_i)_{\alpha}} (\det \msc{Q})^{-\frac{1}{2}} \\
    &= \sqrt{\frac{k}{-N}} \prod_i a_i^{-\frac{1}{2}} L_{\alpha}^{-\frac{N+5}{2}}.
\end{align*}
To obtain~\Cref{lem:vanishing} as a special case of~\Cref{lem:asymptotics}, it remains to check the following facts:
\begin{enumerate}
    \item The constants coming from the Quantum Riemann-Roch theorem match via 
        \[\Delta^{\lambda}(z)^* \mbf{1}_{\pt_{\alpha}} = \Delta^{\pt_{\alpha}}(z) \mbf{1}_{\pt_{\alpha}}; \]
    \item We can identify $q^{-\frac{t_{\alpha}}{z}}(-mzD) \msc{I}_{\alpha}(q,z)$ with $z S^{\lambda}(z)^* \mbf{1}|_{\pt_{\alpha}}$ up to a constant factor.\footnote{Here, note that we are really computing the $I$-function for the total space $\msc{O}(-k)$ on $\P$ via the integral, and this is related to our $I$-function by applying $D$ and multiplying by a prefactor.}
    \item The critical value
        \begin{align*}
            u_{\alpha} &= k L_{\alpha} + \sum_{j=1}^N \zeta_N^j t \log(L_{\alpha} - \zeta_N^j t) \\
            &= \int \frac{NL_{\alpha}^N}{L_{\alpha}^N - t^N} \d L_{\alpha} \\
            &= \int L_{\alpha} \frac{\d q}{q}
        \end{align*}
        differs from $\tau_{\alpha}$ by $t_{\alpha} \log q + \mr{const}$.
\end{enumerate}
Now we may use~\Cref{eqn:rmattwisted} and compare it to~\Cref{eqn:rmatrix}.

\subsection{Polynomiality of $R$-matrix entries}%
\label{sub:Polynomiality of R matrix entries}

Set $I_{33} \coloneqq I_{11}$.
Define a normalized basis\footnote{Here, we follow~\cite[\S 5.2]{nmsp2}.} for $\msc{H}_Z$ and its dual by the formula
\begin{align*}
    \varphi_b &\coloneqq I_0(q') I_{11}(q') \cdots I_{bb}(q')H^b; \\
    \varphi^b &\coloneqq\frac{-t^N}{p_k I_0(q') I_{11}(q') \cdots I_{bb}(q')}H^{3-b}.
\end{align*}
We also define a normalized basis for $H^*(\pt_{\alpha})$ by the formula
\begin{align*}
    \bar{\1}_{\alpha} &\coloneqq L(q')^{-\frac{N+3}{2}} \1_{\alpha}; \\
    \bar{\1}^{\alpha} &\coloneqq L(q')^{\frac{N+3}{2}} \1^{\alpha}.
\end{align*}

Now introduce the following matrix elements for the $R$-matrix:
\begin{align*}
    (R_m)_j^b &\coloneqq (R_m \varphi^b, p^j)^M; \\
    (R_m)_j^{\alpha} &\coloneqq L_{\alpha}^{-j-k} (R_m \bar{\1}^{\alpha}, p^j)^M.
\end{align*}

\begin{lem}\label{lem:level0polynomiality}
    When $m < N-3$, we have $(R_m)_j^b = 0$ if $j \not\equiv b+m \pmod{n}$ and both $(R_m)_{b+m}^b$ and $\frac{Y}{t^N} \cdot (R_m)_{b+N+m}^b$ are in $\msc{R}$.
\end{lem}

\begin{proof}
    First, note that~\Cref{eqn:r1level0} implies that
    \begin{equation}\label{eqn:rlevel0vanishing}
    (R_m)_0^b = \delta_{b,0} \delta_{m,0} 
    \end{equation}
    when $m < N-3$  when restricting to $Z$. We will use the quantum connection to compute the other entries. The quantum connections
    \begin{align*}
        (p+zD) S^M(z)^* &= S^M(z)* \cdot A^M; \\
        (H+zD) S^Z(z)^* &= S^Z(z)* \cdot A^Z 
    \end{align*}
    together with~\Cref{eqn:rmatrix} imply that $R(z)^*|_Z$ satisfies
    \begin{equation}\label{eqn:rlevel0recursion}
     (zD+A^Z) R(z)^*x|_Z = R(z)^* \cdot A^M x |_Z 
    \end{equation}
    for all $x \in \msc{H}$ (simply apply $(p+zD)$ to both sides of~\Cref{eqn:rmatrix} and commute the $S^Z$ to the left of the derivative). This implies that
    \[ (R_m)_j^b = (D+C_b)(R_{m-1})_{j-1}^b + (R_{m-1})_{j-1}^{b-1} - c_j^k q (R_k)_{j-N}^b, \]
    where $C_b = D \log(I_0 I_{11} \cdots I_{bb}) \in \msc{R}$ and $c_j^k = (0, \ldots, 0, c_{1,k}, c_{2,k}, c_{3,k}, c_{2,k})$, where $c_{*,k}$ were defined in~\Cref{lem:quantumconnection}. Here, the $-c_j^kq(R_k)_{j-N}^b$ term comes from the explicit formula for $A^M$, the $(D+C_b)(R_{m-1})_{j-1}^b$ comes from $D (R(z)^*p^j)$ (the correction term $C_b$ comes from differentiating $\varphi^b$), and the $(R_m)_{j-1}^{b-1}$ comes from $A^Q$.

    Now the desired results follow from~\Cref{eqn:rlevel0vanishing} by repeated application of~\Cref{eqn:rlevel0recursion} noting that in~\Cref{eqn:rlevel0recursion}, the quantity $j-b-m \mod{N}$ is preserved.
\end{proof}

\begin{lem}\label{lem:rlevel1polynomiality}
    At the isolated points $\pt_{\alpha}$, the quantity $(R_m)_j^{\alpha}$ is independent of $\alpha$ and
    \[ (R_m)_j^{\alpha} \in \Q[Y]_{m+\floor{\frac{j}{N}}}. \]
\end{lem}

\begin{proof}
    The case of $j=0$ is the second part of~\Cref{prop:rmatlevel1}. We then use the formula
    \[ (\1^{\alpha}, R(z)^* p^j)^M = e^{\frac{\tau_{\alpha}}{z}} \Delta^{\pt_{\alpha}}(z)^* (\1^{\alpha}, S^M(z)^*p^j) \]
    to obtain the recursive formula
    \[ (\1^{\alpha}, R(z)^* p^j) = D_{L_{\alpha}} (\1^{\alpha}, R(z)^* p^{j-1}) - c_j^k q (\1^{\alpha}, R(z)^* p^{j-N}). \]
    Using the derivatives in the first part of~\Cref{rmk:derivatives} (and noting that $X = 1-Y$), we obtain
    \begin{align*}
        (R_m)_j^{\alpha} ={}& \ab(D-\frac{1}{N} \ab(\frac{N+3}{2}-j+m)(1-Y)) (R_{m-1})_{j-1}^{\alpha} \\
        &+ (R_m)_{j-1}^{\alpha} + \frac{c_j}{r}(1-Y)(R_m)_{j-N}^{\alpha}.
    \end{align*}
    Using the base case $j=0$ and the recursive formula (note that $DY = Y(Y-1)$ increases the degree by $1$), we obtain the desired result.
\end{proof}

\begin{lem}
    Define
    \begin{align*}
        V(z,w) &\coloneqq \sum_{\alpha} \frac{e_{\alpha}\otimes e^{\alpha} - R(z)^{-1} e_{\alpha} \otimes R(w)^{-1} e^{\alpha}}{z+w} \\
        &\eqqcolon \sum_{s,t \geq 0} V_{st} z^s w^t.
    \end{align*}
    Then the coefficients $V_{st}$ take the form
    \begin{align*}
        V_{st} ={}& \sum_{a,b=0}^3 ( V_{st} )^{ab} \varphi_a \otimes \varphi_b + \sum_{b=0}^3 \sum_{\alpha=1}^N L_{\alpha}^{2-b-s-t} ( V_{st} )^{\alpha b} \bar{\1}_{\alpha} \otimes \varphi_b \\
        &+ \sum_{\alpha,\beta = 1}^N \sum_j L_{\alpha}^{j-s} L_{\beta}^{2-j-t} ( V_{st} )^{\alpha\beta;j} \bar{\1}_{\alpha} \otimes \bar{\1}_{\beta},
    \end{align*}
    where $\frac{Y}{t^N} ( V_{st} )^{\alpha\beta;j} \in \Q[Y]$ and $\frac{Y}{t^N} ( V_{st} )^{\alpha b}, \frac{Y}{t^N}( V_{st} )^{ab} \in \msc{R}$.  
\end{lem}

\begin{proof}
    The proof is the same as~\cite[Lemma 5.5]{nmsp2}.
\end{proof}

\section{Structure of the GW potential}%
\label{sec:Structure of the GW potential}

\begin{thm}\label{thm:polynomiality}
    For any $g,n$ such that $2g-2+n > 0$, we have
    \[ P_{g,n} \in \msc{R}. \]
\end{thm}

\begin{proof}
    First, by~\Cref{thm:yukawa}, we obtain
    \begin{align*}
        P_{0,3} &= p_k Y I_0^2 I_{11}^3 \ab<\!\ab<H,H,H>\!>_{0,3}^{Z} \\
        &= 1 \in \msc{R}.
    \end{align*}
    We also note that
    \begin{align*}
        D(P_{g,n}) ={}& D\ab(\frac{(p_k Y)^{g-1}I_{11}^n}{I_0^{2g-2}}) \ab(Q\odv{}{Q})^n F_g(Q)\Big\vert_{Q = qe^{\frac{I_1}{I_0}}} \\ 
        &+ \frac{(p_k Y)^{g-1}I_{11}^n}{I_0^{2g-2}} \ab(Q \odv{}{Q}) F_g(Q)\Big\vert_{Q = qe^{\frac{I_1}{I_0}}} \cdot I_{11} \\
        ={}& ( (g-1)(Y-1-2B)+nA )P_{g,n} + P_{g,n+1},
    \end{align*}
    which implies that
    \[ P_{g,n+1} = (D+(g-1)(2B+1-Y)-nA) P_{g,n}. \]
    Because $\msc{R}$ is closed under $D$, we see that if $P_{g,n} \in \msc{R}$, so is $P_{g,n+1}$. We will now fix the second base case $P_{1,1}$.

    Note that in the case of the quintic, $P_{1,1}$ was already known by work of Zinger~\cite{reducedgenus1}. In our case, a direct computation of $\ab<p>_{1,1}^{[0,1]} = \mr{const}$ using~\Cref{lem:level0polynomiality} and~\Cref{lem:rlevel1polynomiality} yields
    \begin{equation}\label{eqn:p11}
     P_{1,1} = -\frac{1}{2} A + \ab(\frac{\chi(Z)}{24} - 2)B - \frac{1}{12}X + a_{1,k}, 
    \end{equation}
    where $\chi(Z_6) = -204$, $\chi(Z_8) = -296$, $\chi(Z_{10}) = -288$, and $a_{1,k} = \int_{Z_k} c_2(Z_k) \cup H$ is given in~\Cref{tab:label}. The stable graphs entering this computation are given in~\Cref{fig:03graphloop}. For example, at level $0$, the graph with one genus one vertex and one leg contributes
    \[ P_{1,1} - \frac{\chi(Z)}{24}B \]
    and the graph with one genus zero vertex with one leg and one loop contributes\footnote{Our $r_{0,k}$ is related to the $r_0$ of~\cite{hkq}, so we adopt their notation.}
    \[ \frac{1}{2} A + 2B + r_{0,k} X, \]
    where $r_{0,6} = \frac{13}{36}$, $r_{0,8} = \frac{11}{32}$, and $r_{0,10} = \frac{3}{10}$.
    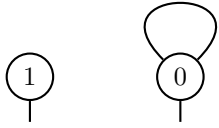
\begin{figure}[htpb]
    \begin{center}
    \begin{tikzpicture}[scale=1, transform shape, every loop/.style={}]
    \node[circle,draw=black,thick] (B) at (-1,0) {$1$};
    \draw[-,thick] (B) -- (-1,-0.6);
    \node[circle,draw=black,thick] (A) at (1,0) {$0$};
    \draw[-,thick] (A) -- (1,-0.6);
    \path[every node/.style={font=\sffamily\small},thick]
        (A)   edge[loop] (A);
    \end{tikzpicture}
    \end{center}
    \caption{Stable graphs for $g=1, n=1$}%
    \label{fig:03graphloop}
    \end{figure}

    Fix a genus $g \geq 2$. Assume that for any $h < g$ and any $n$ such that $2g-2+n > 0$, we have $P_{h,n} \in \msc{R}$. We will consider the ``normalized'' genus $g$ MSP generating function
    \[ (p_k Y )^{g-1} \ab<\ >_{g,0}^{[0,1]} \]
    when $g \geq 2$. We have $(p_k Y)^{g-1} \ab<\ >_{g,0}^{[0,1]} \in \Q[Y]$ by~\Cref{thm:polynomiality01} and the fact that $q = \frac{1}{r}(1-Y^{-1})$.

    We now consider the graph sum formula of~\Cref{thm:01cohft}. The leading graph is a single genus $g$ vertex at level $0$ with no marked points, which contributes
    \[ \frac{(p_k Y)^{g-1}}{I_0^{2g-2}} F_g = P_g. \]
    We consider all of the remaining graphs. Note that
    \[ \sum_v (g_v-1) + \ab|E| = g-1, \]
    so we can distribute the $Y^{g-1}$ factor into the vertices and edges. We will prove polynomiality of the contribution of each such graph.
    \begin{itemize}
        \item At each edge, the contribution has the form $Y \cdot ( V_{st} )^{ab}$, $Y \cdot ( V_{st} )^{\alpha b}$, or $Y \cdot ( V_{st} )^{\alpha\beta; j}$, which are independent of the hour (in the last two cases) and lie in $\msc{R}$;
        \item At each level $0$ vertex, the contribution is 
            \[ Y^{g_v-1} \int_{\ol{\msc{M}}_{g_v, n_v}} \ab[\varphi_{a_1}\bar{\psi}_1^{m_1}, \ldots, \varphi_{a_{n_v}} \bar{\psi}_{n_v}^{m_{n_v}}]_{g_v, n_v}^{Z, T}. \]
            This translated correlator vanishes unless $\sum a_i + m_i = n_v$ for virtual dimension reasons. Using the string and dilaton equations (and dimension constraints), this reduces to a constant multiple of $P_{g_v, m}$ (here, it is helpful to note that ancestors on $\ol{\msc{M}}_{0,3} = \pt$ vanish). Because $g_v < g$, this contribution lies in $\msc{R}$;

        \item At each level $1$ vertex, the contribution is\footnote{Note that the $L^{n\frac{N+3}{2}}$ in the statement of~\Cref{thm:01cohft} is absorbed by using the basis $\bar{\mbf{1}}_{\alpha}$ instead of $\mbf{1}_{\alpha}$.}
            \begin{align*}
                &Y^{g_v-1} \int_{\ol{\msc{M}}_{g_v, n_v}} \ab[L_{\alpha}^{j_1-m_1} \bar{\psi}_1^{m_1}, \ldots, L_{\alpha}^{j_{n_v} - m_{n_v}} \bar{\psi}_{n_v}^{m_{n_v}}]_{g_v, n_v}^{\pt_{\alpha}, T} \\
                &= \sum_s \frac{L^{\frac{3}{2}(2g_v-2)}}{s!} \ab<L_{\alpha}^{j_1-m_1}\bar{\psi}_1^{m_1}, \ldots, L_{\alpha}^{j_{n_v} - m_{n_v}} \bar{\psi}_{n_v}^{m_{n_v}}, \prod_{a=1}^s \wt{T}_{\alpha}(\bar{\psi}_{n_v+a})>_{g_v, n_v+s}.
            \end{align*}
            After summing over $\alpha$, this contribution will lie in $\msc{R}$. 
            \begin{enumerate}
                \item If we write
                    \[ \wt{T}_{\alpha}(\bar{\psi}_{n_v+a}) = \sum_{\ell_a} (\wt{T}_{\alpha})_{\ell_a} \bar{\psi}_{n_v+s}^{a_s+1}, \]
                    then the correlator corresponding to a given choice of monomials is nonzero if and only if
                    \[ \sum_{i=1}^{n_v} m_i + \sum_{a=1}^s (\ell_a + 1) = 3g_v-3+n_v+s. \]
                \item Using the fact that $L_{\alpha}^{\ell} (\wt{T}_{\alpha})_{\ell} = (R_{\ell})_j^{\alpha} \in \Q[Y]$ by~\Cref{lem:rlevel1polynomiality}, the total factor involving $L_{\alpha}$ is $L_{\alpha}^{\ab(\sum_i j_i) - n_v}$, which becomes a multiplicative factor of the contribution after summing over $\alpha$.

                    We now assume that $N$ is prime. After summing over all $\alpha$, this is nonzero only if the resulting sum has no roots of unity, or in other words,
                    \[ \sum_{i=1}^{n_v} j_i \equiv n_v \pmod{N}. \]
                    Using the fact that for each edge, one end contributes $L_{\alpha}^j$ and the other end contributes $L_{\beta}^{2-j}$, the total factors in the graph become
                    \[ \prod_{v} L_{\alpha(v)}^{\sum_{i=1}^{n_v} j_i - n_v} = \prod_{v} Y^{\frac{\sum_{i=1}^{n_v} j_i - n_v}{N}} = 1. \]
            \end{enumerate}
            We now see that all contributions from graphs besides the leading one lie in $\msc{R}$. Therefore we have
            \[ P_g = \sum_{\Gamma \text{ non-leading}} \frac{1}{\ab|\Aut\Gamma|} \on{Cont}_{\Gamma} \in \msc{R}. \qedhere \]
    \end{itemize}
\end{proof}

\begin{thm}\label{thm:analytic}
    For any $g \geq 0$, the function 
    \[ \mc{F}_g(Q) \coloneqq \begin{cases}
        \ab(Q \odv{}{Q})^3 F_0(Q) & g=0 \\
        \ab(Q \odv{}{Q}) F_1(Q) & g=1 \\
        F_g(Q) & g\geq 2 
    \end{cases}
    \]
    is an analytic function in $Q$ near $Q = 0$. 
\end{thm}

\begin{proof}
    Note that
    \[ P_g = \frac{(p_k Y)^{g-1}}{I_0^{2g-2}} F_g \in \msc{R}. \]
    Then because $I_0, I_1, I_2, I_3$ are all analytic funtions in $q$ with radius of convergence $\frac{1}{r}$, so are all of the generators $A_m$ and $B_m$. So is the mirror map $\frac{I_1}{I_0}$ (here note that the constant term of $I_0$ is $1$). Because the map $q \mapsto qe^{\frac{I_1}{I_0}}$ is a biholomorphism between neighborhoods of $0$ (here, use the fact that $I_1(0) = 0$, so the mirror map is $q(1+\cdots)$), we see that all of the components of the $I$-function, the mirror map, and all elements of $\msc{R}$ are analytic functions of $Q$ near $Q=0$. We conclude that
    \[ \mc{F}_g(Q) = \begin{cases}
        p_k Y I_{11}^{-3} I_0^{-2} P_{0,3} & g=0 \\
        I_{11}^{-1} P_{1,1} & g=1 \\
        p_k^{1-g}I_{11}(1-rq)^{g-1}I_0^{2g-2} & g \geq 2
    \end{cases}
    \]
    is an analytic function of $Q$ near $Q = 0$.
\end{proof}

\appendix

\section{The normalized Yukawa coupling}
\label{sec:yukawa coupling}

We will now prove that the normalized Yukawa coupling is a rational function using the same strategy as in the proof of~\cite[Theorem 2]{zagierzinger}. For simplicity, we will only give the proof for $k=6$. The proofs for $k=8, 10$ are essentially identical. In the rest of this section, we will use the notation of Zagier-Zinger. In particular, their $I_p$ is our $I_{pp}$.

Let $D_w \coloneqq w + x \odv{}{x}$. Then, following Zagier-Zinger, define
\[ \mbf{M}F(w,x) \coloneqq w^{-1} D_w\ab(\frac{F(w,x)}{F(0,x)}) \]
whenever $F \in 1+x\Q(w)\llbracket x \rrbracket$ has coefficients holomorphic at $w=0$.

\begin{thm}
  Let
  \[ \msc{F}(w,x) \coloneqq \sum_{d=0}^{\infty} x^d \frac{\prod_{r=1}^{6d} (6w+r)}{2^d\prod_{r=1}^d ((w+r)^5(2w+2r-1))} \]
  and define $\msc{F}_p(w,x) \coloneqq \mbf{M}^p \msc{F}(w,x)$ and $I_p(x) = \msc{F}_p(0,x)$. Then
  \begin{enumerate}
    \item $I_0(x) \cdots I_4(x) = Y$;
    \item For $0 \leq p \leq 4$, we have $I_p = I_{4-p}$.
  \end{enumerate}
\end{thm}

\begin{proof}
  First, note that $\msc{F}$ satisfies the Picard-Fuchs equation
  \[ D_w^4 - \frac{1}{2^2} 6^2 x \prod_{j=1,2,4,5}(6 D_w + j) = w^4. \]
  Applying~\cite[Corollary 1]{zagierzinger}, we obtain
  \[ \sum_{s=0}^{4-p} C_s^{(p)}(x) D_w^s \msc{F}_p(w,x) = w^{4-p}, \]
  where the coefficients $C_s^{(p)}(x)$ can be computed recursively. The top one is given by
  \[ C_{4-p}^p = Y^{-1} I_0(x) \cdots I_{p-1}(x), \]
  so setting $p=4$, we obtain
  \[ Y^{-1} I_0(x) \cdots I_3(x) \msc{F}_4(w,x) = 1. \]
  Therefore, $\msc{F}_4(w,x) = I_4(x)$ is independent of $w$ and
  \[ I_0(x) \cdots I_4(x) = Y. \]
  This proves~\Cref{thm:yukawa}.

  We can now repeatedly solve ODEs to reconstruct $\msc{F}_p(w,x)$ from their initial values $I_p(x)$. In particular, we obtain
  \[ w^{-4}\msc{F}(w,x) = I_0 D_w^{-1}I_1 D_w^{-1}I_2 D_w^{-1}I_3 D_w^{-1}I_4. \]
  Matching coefficients of $x^d$, we obtain
  \[\frac{6^{-1} 2^{-d} \prod_{r=0}^{6d}(6w+r)}{\prod_{r=0}^d (w+r)^5 \prod_{r=1}^d(2w+2r-1)} = \sum_{\substack{d_0,\ldots,d_{n-1} \geq 0 \\ d_0 + \cdots + d_{n-1}=d}} \frac{c_0(d_0) \cdots c_4(d_{n-1})}{(w+d_1 + \cdots + d_{4})\cdots(w+d_{4})},\]
  where $c_p(d)$ is defined by
  \[ I_p(x) = \sum_d c_p(d) x^d. \]
  Considering the leading terms of the right hand side, meaning the permutations of $(d,0,\ldots,0)$, we obtain
  \begin{align*}
    \sum_{p=0}^4 \frac{c_p(d)}{w^{4-p}(w+d)^p} =&{}\frac{6^{-1} 2^{-d} \prod_{r=0}^{6d}(6w+r)}{\prod_{r=0}^d (w+r)^5 \prod_{r=1}^d(2w+2r-1)} \\
    &- \sum_{\substack{d_0,\ldots,d_{n-1} > 0 \\ d_0 + \cdots + d_{n-1}=d}} \frac{c_0(d_0) \cdots c_4(d_{n-1})}{(w+d_1 + \cdots + d_{4})\cdots(w+d_{4})}.
  \end{align*}
  We will now prove $c_p(d) = c_{4-p}(d)$ by induction. The base case $d=0$ is satisfied because $c_p(0) = I_p(0) = 1$. Now suppose that $c_p(d') = c_{4-p}(d')$ whenever $d' < d$. Then the right hand side is invariant under $w \mapsto -w-d$. For the first term, we simply see that both numerator and denominator have an odd number of terms. For the second term, we see the invariance by reordering $d_p \mapsto d_{4-p}$ and noting the number of terms in the denominator is even. Therefore, the left hand side is also invariant under $w \mapsto -w-d$, so we must have $c_p(d) = c_{4-p}(d)$, as desired. Reforming the generating series, we obtain $I_p(x) = I_{4-p}(x)$.
\end{proof}

\begin{rmk}
    Our $\msc{F}$ is Zagier-Zinger's $\wt{\msc{F}}$. If we instead consider 
    \[ \msc{F}'(w,x) \coloneqq \sum_{d=0}^{\infty} x^d \frac{\prod_{r=1}^{6d} (6w+r)}{\prod_{r=1}^d (2(w+r)^5(2w+2r-1) - 4w^6)}, \]
    we obtain the result that the $\msc{F}'_p(w,x)$ are periodic for $M$ with period $6$ (or $k$ in general) using the same method as in the proof of~\cite[Theorem 2]{zagierzinger}.
\end{rmk}

\printbibliography

\end{document}